\documentclass[a4paper,12pt]{article}

\usepackage[affil-it]{authblk}

\usepackage{amsmath,amssymb}
\usepackage{amsthm}
\usepackage{mathtools}
\usepackage{graphicx}

\usepackage{csquotes} %\enquote{text}

\usepackage[]{algorithm2e}
\usepackage{enumitem}

\usepackage{tcolorbox}

\newcommand{\dd}{\mathrm{d}}

\usepackage{hyperref}

\usepackage[margin=2.5cm]{geometry}

\graphicspath{ {./fig/} }

\usepackage[font=small, format=hang, labelfont=bf]{caption}
\usepackage[font=footnotesize, format=hang, margin={0.15\textwidth,0.1\textwidth}]{subcaption}
\usepackage{float}

\usepackage{booktabs}

\newcommand{\sca}{0.36}
\renewcommand{\d}{\mathrm{d}}

\newtheorem{theorem}{Theorem}[section]

\theoremstyle{definition}

\theoremstyle{remark}
\newtheorem{remark}[theorem]{Remark}
\newtheorem{remark*}{Remark}

\numberwithin{equation}{section}

\title{Numerical study of Bose--Einstein condensation in the Kaniadakis--Quarati model for bosons}

\author[ ]{ J.\;A.\;Carrillo\,\thanks{\tt carrillo@imperial.ac.uk}\hspace{-3pt}}
\author[ ]{K.\;Hopf\,\thanks{\texttt{hopf@wias-berlin.de} (corresponding author)}\hspace{-3pt}}
\author[ ]{M.-T.\;Wolfram\,\thanks{\tt m.wolfram@warwick.ac.uk}}

\affil[$*$]{\small Department of Mathematics, Imperial College London,
  
       London SW7 2AZ, UK \vspace{.75em}}
\affil[$\dag$,$\ddag$]{Mathematics Institute, University of Warwick,
  
  Coventry CV4 7AL, UK}

\date{\small\today}
\makeindex
\date{}

\begin{document}
 \maketitle
 
 \begin{abstract}\noindent
Kaniadakis and Quarati (1994) proposed a Fokker--Planck equation with quadratic drift as a PDE model for the dynamics of bosons in the spatially homogeneous setting. It is an open question whether this equation has solutions exhibiting condensates in finite time. The main analytical challenge lies in the continuation of exploding solutions beyond their first blow-up time while having a linear diffusion term. We present a thoroughly validated time-implicit numerical scheme capable of simulating solutions for arbitrarily long time, and thus enabling a numerical study of the  condensation process in the Kaniadakis--Quarati model. We show strong numerical evidence that above the critical mass rotationally symmetric solutions of the Kaniadakis--Quarati model in $3$D form a condensate in finite time and converge in entropy to the unique minimiser of the natural entropy functional. 
Our simulations further indicate that the spatial blow-up profile near the origin follows a universal power law and that transient condensates can occur for sufficiently concentrated initial data. 
  
  \bigskip
  
  \textbf{Keywords}: nonlinear Fokker--Planck equation; Bose--Einstein condensation; 
  entropy decay; implicit-in-time discretization; Lagrangian scheme.
 \end{abstract}
  
 \section{Introduction}\label{sec:introSim}

 In this paper we consider the following family of nonlinear Fokker--Planck equations
 \begin{align}\label{eq:befp}
   \partial_t f & =   \Delta_v  f+\mathrm{div}_v(v f(1+ f^\gamma)),\quad t>0,\;v\in\mathbb{R}^d,
   \\\nonumber f(0,\cdot) & =  f_0\ge0,
 \end{align}
 where $\gamma>0$ is a fixed parameter and $f=f(t,v)\ge0$.
 We are particularly interested in the case $\gamma=1$,
 in which equation~\eqref{eq:befp} is known as the \textit{Kaniadakis--Quarati model} for bosons (KQ).
It was introduced by Kaniadakis and Quarati~\cite{kaniadakis_classical_1994} as a model for quantum particles following Bose--Einstein statistics, obtained by adapting accordingly the transition probability rates in the kinetic model.

\paragraph{Physical background (\texorpdfstring{$\gamma=1$}{}).} 
 
The feature in which KQ differs from the linear Fokker--Planck equation consists in the additional factor $(1+f)$ in the drift term. This factor, leading to a nonlinear equation, arises from the assumption of 
indistinguishability of identical quantum particles. Indeed, in contrast to classical mechanics, in a quantum system of identical and indistinguishable particles, the presence of particles in a given energy state influences the probability of further quantum particles joining that state. Here we are interested in systems of bosons, whose wave function is symmetric with respect to permutations of particles. This results in an increase in the transition probability, which is encoded, in the continuum model, in the extra factor $(1+f)$. For KQ the choice $d=3$ is the physically most interesting space dimension. In this case the problem exhibits a finite critical mass $m_c$ above which condensates are expected to emerge in finite time, see below for more details. However, in the literature little is known about the possible formation of condensates in 3D KQ.

\paragraph{Variational structure and steady states.}
Equation~\eqref{eq:befp} has a natural \textit{entropy functional}, given by
\begin{align*}
 \mathcal{H}(f):=\int \left(\frac{|v|^2}{2}f+\Phi(f)\right)\,\d v,
\end{align*}
where $\Phi(f):=\frac{1}{\gamma}\int_0^f\log\left(\frac{s^\gamma}{1+s^\gamma}\right)\d s$ and thus $\Phi''(f)=1/h(f)$ for $h(s):=s(1+s^\gamma)$.
Indeed, formally, equation~\eqref{eq:befp} can be rewritten as
\begin{align}\label{eq:gradflow}
 \partial_tf = \nabla\cdot\left(h(f)\nabla\frac{\delta \mathcal{H}}{\delta f}(f)\right),  
\end{align}
where $\frac{\delta \mathcal{H}}{\delta f}$ denotes the variational derivative of $\mathcal{H}$.
Thus, for any sufficiently regular, positive (and hence mass conserving) solution $f=f(t,v)$ of eq.~\eqref{eq:befp}, one obtains the \textit{entropy dissipation identity}
 \begin{align}\label{eq:ediss}
   \frac{\d}{\d t}\mathcal{H}(f)=-\int h(f)\left|\nabla\frac{\delta\mathcal{H}}{\delta f}(f)\right|^2\,\d v.
 \end{align}
Notice, however, that due to the presence of the (quantum correction) term $s^{\gamma}$ in the definition of $h(s)$, equation~\eqref{eq:gradflow} is not a gradient flow of the functional $\mathcal{H}$ with respect to the classical Wasserstein metric. 
The mobility $h(s)$ associated to the nonlinear continuity equation~\eqref{eq:gradflow}
is convex leading to well-known issues of ill-defined Wasserstein-like metrics to render rigorous the gradient flow structure~\cite{DNS} in contrast to the Fermi--Dirac case \cite{CLR,CLSS}. 

We observe that, given a sufficiently regular positive function $f$, the RHS of equation~\eqref{eq:ediss} is strictly negative unless $\nabla\frac{\delta\mathcal{H}}{\delta f}(f)=0$. The regular solutions of this equation are henceforth referred to as the \textit{steady states} associated with problem~\eqref{eq:befp}. They are explicitly given by 
\begin{align}\label{eq:ss}
  f_{\infty,\theta}(v)=\left(\mathrm{e}^{\gamma(\frac{|v|^2}{2}+\theta)}-1\right)^{-1/\gamma},\quad\theta\ge0.
\end{align}
Notice that $f_{\infty,\theta}$ is smooth and integrable for $\theta>0$, and the family $\{f_{\infty,\theta}\}$ is strictly ordered and approaches $f_c:=f_{\infty,0}$ from below as $\theta\searrow0$. Furthermore, letting $m_c:=\int f_c$, the map $(0,\infty)\ni\theta\mapsto m_\theta:=\int f_{\infty,\theta}\in(0,m_c)$ is a bijection, and $m_c<\infty$ if and only if $\gamma>\frac{2}{d},$ i.e.\;if and only if the problem is $L^1$-supercritical. While $f_{\infty,\theta}$ is the unique minimiser of $\mathcal{H}$ among non-negative integrable functions of mass $m=m_\theta$, for $m>m_c$ the problem of minimising $\mathcal{H}$ under mass constraint does not have a regular solution. 
Since $\Phi$ is sublinear at infinity, the natural extension  $\mathcal{\widetilde H}$ of the entropy functional to the set of finite non-negative Borel measures $\mathcal{M}^+_b$ is given by
$$
\mathcal{\widetilde H}:\quad\mu\mapsto\int\left(\frac{|v|^2}{2}\mu(\d v)+\Phi(f)\,\d v\right),
$$
where $f$ denotes the density of the absolutely continuous part of $\mu$. The extension
 is convex and lower-semicontinuous with respect to {weak-star} convergence in $\mathcal{M}$~\cite{demengel_convex_1984,abdallah_minimisation_2011}. In~\cite{abdallah_minimisation_2011} it is shown via an explicit calculation that the extended functional has a  unique minimiser among finite non-negative measures of mass $m>m_c$, which is given by 
\begin{align*}
  f_c\cdot\mathcal{L}^d+(m-m_c)\delta_0.
\end{align*}
The above comments on the entropy functional and the steady states of equation~\eqref{eq:befp} apply to the problem posed on the whole space $\mathbb{R}^d$ (assuming sufficient decay as $|v|\to\infty$) as well as to the problem on a centred ball $B(0,R_1)$ subject to no-flux boundary conditions.
 
 \paragraph{Dynamics of the Kaniadakis--Quarati model.}
As noted in~\cite{carrillo_finite-time_2019}, in the $L^1$-subcritical case, $d=1$,  KQ is globally wellposed in the classical sense for sufficiently regular initial data, and solutions converge to equilibrium at an exponential rate~\cite{carrillo_1d_2008}. In the $L^1$-critical case, $d=2$, solutions are also globally regular and converge to equilibrium~\cite{canizo_fokkerplanck_2016}---with an exponential rate in the spatially isotropic case $f(t,v)=g(t,|v|)$. The approach in~\cite{canizo_fokkerplanck_2016} exploits the fact that $2$D KQ in isotropic coordinates can be transformed to a linear Fokker--Planck equation, which leads to explicit solutions also for the nonlinear equation.
For $3$D KQ Toscani~\cite{toscani_finite_2012} proved via contradiction the existence of solutions blowing up in finite time. Finite-time blow-up in this reference is obtained for any solution of sufficiently large mass $m$ (above a technical threshold far larger than the critical mass), but also for solutions of arbitrarily small mass provided they are initially sufficiently concentrated near the origin. 
Formal results on the dynamics of isotropic solutions to 3D KQ based on matched asymptotic expansions have been obtained in~\cite{sopik_dynamics_2006}. Our numerical simulations will qualitatively confirm some of the main findings in~\cite{sopik_dynamics_2006}, which suggests that the dynamics depicted in this reference give a good hint at the typical behaviour of solutions.
Our numerical experiments will, however, also indicate that the dynamics may, in general, display a richer variety of phenomena. The formal considerations in~\cite{sopik_dynamics_2006} rely on the assumption of a sufficiently spread out initial datum $f_0$. We would also like to emphasize that in contrast to~\cite{sopik_dynamics_2006} our scheme allows for a numerical study beyond the first blow-up time.

\paragraph{$L^1$-supercritical Fokker--Planck model for bosons in 1D.} The one-dimensional case of equation~\eqref{eq:befp} with $\gamma>2$ was recently studied in the ref.~\cite{carrillo_finite-time_2019}, both on the entire line as well as on a centred interval subject to zero-flux boundary conditions.
	We would like to point out that the successful numerical experiments reported in the present manuscript,
	which are based on the equation for the pseudo-inverse cumulative distribution function
	\begin{align}\label{eq:defPsI}
	u(x) = \inf\left\{r:\int_{\{r'\le r\}}f(r')\,\dd r'\ge x\right\},\quad x\in(0,\|f\|_{L^1}),
	\end{align}
	of the original density $f=f(t,\cdot)$ (cf.~equation~\eqref{eq:equ} in Section~\ref{sec:icdf} below)
	triggered the rigorous analysis in~\cite{carrillo_finite-time_2019}, which is itself based on this reformulation. 
	Let us briefly review those results of~\cite{carrillo_finite-time_2019} which are of relevance for the present paper: the authors obtain global-in-time existence and uniqueness of solutions $u$ in the viscosity sense for initial data corresponding to sufficiently regular positive densities $f_0$ of finite mass $m$. These solutions are non-decreasing in the mass variable, here denoted by $x$. It is further shown that such solutions $u=u(t,x)$ are smooth away from $\{u=0\}$ and that 
	the push-forward measure $u(t,\cdot)_\#\mathcal{L}^1_{|(0,m)}=:\mu(t)\in\mathcal{M}^+_b$, generalising the problem in the original variables, has the the form
\begin{align*}
 \mu(t) =  f(t,\cdot)\cdot\mathcal{L}^1+ x_p(t)\delta_0,
\end{align*}
where the map $t\mapsto x_p(t):=\mathcal{L}^1(\{u(t,\cdot)=0\})$ is continuous and the function $f(t,\cdot)\in L^1_+$ is smooth away from the origin, where it satisfies equation~\eqref{eq:befp} in the pointwise sense.
Moreover, whenever the density $f(t,\cdot)$ is unbounded at the origin, its spatial blow-up profile has the form
\begin{align}\label{eq:profile}
  f(t,v)=f_c(v)\cdot(1+O(|v|)=c_\gamma|v|^{-\frac{2}{\gamma}}\left(1+O(|v|)\right)\quad\text{as }|v|\to0,\qquad c_\gamma=\left(2/\gamma\right)^\frac{1}{\gamma}.
\end{align}
See~\cite{hopf_thesis} for a refinement of~\eqref{eq:profile}. The above framework makes it possible to extend entropy methods globally in time and to deduce convergence to the measure of the same mass which minimises the entropy. In the case $m>m_c$, the minimiser has a positive Dirac mass at the origin, and the solution must eventually have a non-trivial condensate component.
On the other hand, if $m<m_c$, the minimiser is smooth, and from the bound~\eqref{eq:profile} it can easily be deduced that in this case there exists $T\in(0,\infty)$ such that $x_p(t)=0$ for all $t\ge T$ (see~\cite[Cor.~3.16]{carrillo_finite-time_2019}).
From this observation combined with an adaptation of the finite-time blow-up argument in~\cite{toscani_finite_2012} one infers the existence of solutions whose condensate component $x_p=x_p(t)$ is not identically zero but compactly supported in $(0,\infty)$ (see~\cite[Cor.~3.18]{carrillo_finite-time_2019}). We refer to this phenomenon as a \textit{transient condensate}. Below we will see that the $L^1$-supercritical case in 1D of the family of nonlinear Fokker-Planck equations~\eqref{eq:befp}, corresponding to~$\gamma>2$, appears to be a good 
caricature for the dynamical behaviour of the physically interesting case of the 3D KQ model in radial coordinates.

Let us finally mention that equation~\eqref{eq:befp} in 1D and without the diffusion term was analysed in \cite{CDT} showing that condensates always form in finite time and that their mass is increasing in time so that, once formed, they never dissolve. The results reported here and in \cite{carrillo_finite-time_2019} show the genuine countereffect of linear diffusion on condensation leading to transient condensates and non-monotonic behaviour of the condensate part $x_p(t)$, proved in one dimension for $\gamma>2$ and conjectured in the three dimensional case for $\gamma=1$.

\paragraph{Main numerical findings.}
The main purpose of this work is to provide strong numerical evidence for the existence of solutions to 3D KQ forming a Bose--Einstein condensate in finite time. Our numerical results suggest that any rotationally symmetric solution above the critical mass will eventually have a non-trivial condensate component. From our simulations a rather clear picture of the dynamical properties of KQ in $3$D in the isotropic case will emerge: the long-time asymptotics will be identified, which the numerical solution converges to in entropy at an exponential rate. Numerical evidence is provided for the possibility of  the condensed part failing to be monotonic in time and for even dissolving completely. 
Before investigating KQ in $3$D, we will apply the numerical scheme to the caricature of the $L^1$-supercritical case in 1D, i.e.\;\eqref{eq:befp} with $\gamma>2$, in order to numerically reproduce the analytical results established rigorously in~\cite{carrillo_finite-time_2019}, see Section~\ref{ssec:sim1D}. Since non-stationary explicit solutions are not available in $1$D, the $1$D scheme (in the $L^1$-supercritical case) will be validated by numerically analysing the convergence behaviour under mesh refinement with respect to a reference solution on a very fine mesh. 
Concerning the scheme for rotationally symmetric solutions of KQ we perform a validation in 2D, where explicit solutions are available.

 \paragraph{Numerical scheme.} The proposed numerical scheme is based on the variational formulation of equation~\eqref{eq:befp} using a mass transportation Lagrangian approach.    It  is motivated by the approach in~\cite{blanchet_convergence_2008,carrillo_numerical_2016}, where the gradient flow with respect to the Wasserstein distance is expressed in terms of the inverse of the cumulative distribution functions. Inherent in this approach is the conservation of mass property, which follows by construction. 
We would like to emphasize that concerning the Kaniadakis--Quarati model in 3D studied in the present work, far less is known rigorously as compared to the equations simulated in~\cite{blanchet_convergence_2008,carrillo_numerical_2016} (porous medium equation, critical Keller--Segel) which have been exhaustively studied in the literature. This is partially explained by the fact that the variational structure for this problem cannot directly be exploited by resorting to established tools from optimal transportation theory.
In fact, the potential difficulty in our situation lies in the circumstance that we do not have the Wasserstein gradient flow structure in a rigorous sense. We will, however, see that this precise structure is not required and our proposed scheme will be shown to preserve in particular the entropy decay property (rigorously in 1D and 2D for the semidiscrete case, see~Section~\ref{ssec:semidiscrete}). Our numerical scheme is able to go beyond the first blow-up time and allows for exploring the qualitative behaviour after blow-up: blow-up profile, transient condensates and entropy decay. These good numerical properties, consistent in 1D with the existing theory,  reassure us in our numerical findings in Section~\ref{ssec:sim3D} concerning the 3D isotropic case. The fact that our numerical experiments clearly support the conjecture that the qualitative behaviour of condensates in 1D proven in \cite{carrillo_finite-time_2019} is expected in the most realistic case of 3D radially symmetric initial data can be regarded as  the main contribution of this paper.
\color{black}

There has been an increased interest in related structure preserving Lagrangian schemes in the last years, see for example~\cite{gosse_toscani_2006,blanchet_convergence_2008, carrillo_moll_2009, matthes_osberger_2014, CB16, carrillo_numerical_2016,carrillo2017blob,carrillo_lagrangian_2018}. The numerical analysis of these schemes is still underdeveloped with partial results in~\cite{blanchet_convergence_2008,matthes_osberger_2014,CG,carrillo2017blob,carrillo_lagrangian_2018}. 
Let us finally point out that free energy decaying numerical schemes in the original variables based on finite volume schemes have been proposed in~\cite{CCH15,PZ18,ABPP,BCH} and references therein. These schemes fail to go beyond the blow-up time since they cannot resolve the presence of Dirac concentrations while  accurately following the evolution of the smooth part of the solution.

 \paragraph{Comparison with other models for Bose--Einstein condensation.}
 
There are many other models in the literature which have been suggested in the context of Bose--Einstein condensation. Of particular interest (due to similar phenomena) is a certain class of kinetic equations generally referred to as quantum Boltzmann equations, which, in contrast to classical Boltzmann equations, are derived using Bose--Einstein statistics.
Let us note that for $\gamma=1$ the steady states~\eqref{eq:ss} coincide with the classical Bose--Einstein distributions and the functional $\int\Phi(f)\,\d v$ agrees (up to a sign convention) with the entropy associated to the homogeneous Boltzmann--Nordheim equation for bosons, see~\cite{escobedo_finite_2015,huang_statistical_1963}. In contrast to equation~\eqref{eq:befp}, the Boltzmann--Nordheim equation formally preserves the kinetic energy $\int\frac{|v|^2}{2}f\,\d v$.
 In the last two decades, significant progress has been made in the analysis of the Boltzmann--Nordheim equation in the homogeneous and velocity isotropic case~\cite{escobedo_quantum_2001, escobedo_asymptotic_2004,
lu_boltzmann_2013, escobedo_finite_2015, escobedo_blow_2014, bandyopadhyay_blow-up_2015, lu_2016, lu_strong_2018}. 
To roughly summarise the main results, the authors of the cited references are able to establish the existence of generalised mass- and energy-conserving solutions, which form a Bose--Einstein condensate in finite time and converge, in some sense and under certain conditions, to the entropy minimiser in the large-time limit.
The results in the present paper suggest that in the isotropic case the dynamics of condensation in 3D KQ is in some aspects similar to the one of the Boltzmann--Nordheim equation as described rigorously in the references~\cite{escobedo_finite_2015, escobedo_blow_2014, bandyopadhyay_blow-up_2015,lu_2016, lu_strong_2018}. 
We note that, regarding the nature of singularities, 
in the Boltzmann--Nordheim equation many questions are still open.

Numerical schemes to approximate the Boltzmann--Nordheim equation or quantum Boltzmann equation for bosons have also been devised and used to understand their qualitative properties, see \cite{MP,BPM,HLP,FHJ} and the references therein. However, only few numerical studies attempt to go beyond the first blow-up time (where the velocity distribution ceases to be bounded).  In \cite{semikoz_condensation_1997,semikoz_kinetics_1995,lacaze_dynamical_2001,spohn_kinetics_2010} the authors observe that at the first blow-up time the solution has an integrable power law singularity near the lowest energy state and, in general, there will be a non-trivial flux of particles entering that state. The hypothesis of mass conservation then leads to a law for the time evolution of the condensate component, resulting in a coupled system. The methods do not appear to allow to track in a precise way the evolution after blow-up. Our approach is very different as it does not require  distinguishing between the times where the velocity distribution is bounded and the times where it is unbounded, and enables a detailed study of the dynamics of singular solutions. Let us finally mention that other descriptions have been used both analytically and numerically to study the behaviour beyond condensation in the quantum Boltzmann equation. In some of them the kinetic equation is coupled to a nonlinear Schr\"odinger equation (cubic, Gross--Pitaevski) modelling the evolution of the condensate, see  \cite{ST,BC,B,escobedo_turbulence_2015} and the references therein for further details.

 \paragraph{Plan of the manuscript.}
 The remaining part of this manuscript is structured as follows: in Section~\ref{sec:icdf} we discuss the numerical scheme for the 1D caricature of the 3D KQ given by the $L^1$-supercritical 1D Fokker-Planck equation \eqref{eq:befp} with $\gamma>2$ and its  generalisation to the radial case in higher-dimensions with particular focus on the KQ model, $\gamma=1$. We also briefly discuss the anisotropic case.
  Section~\ref{ssec:sim1D} shows that the proposed numerical scheme does capture the main behaviour after blow-up in the 1D case: condensation, transient condensates for subcritical initial mass and convergence towards equilibrium.  Section~\ref{ssec:val_2D} validates the discretisation of the radial case by comparing to the explicit solutions given in~\cite{canizo_fokkerplanck_2016}. 
In  Section~\ref{ssec:sim3D} we present the simulations of 3D~KQ, which 
allow us to conclude that the caricature given by the $L^1$-supercritical  Fokker-Planck equation~\eqref{eq:befp} in 1D is essentially numerically correct for the 3D KQ model for radial initial data.

\section{Numerical method}\label{sec:icdf}

Since we want our scheme to be able to deal with Dirac masses at the origin, our simulations are not based on the formulation~\eqref{eq:befp}; instead we follow and generalise the ansatz in the ref.~\cite{carrillo_finite-time_2019} considering the equation satisfied by the (pseudo-) inverse cumulative distribution function (cdf) of $f(t,\cdot)$. In higher dimensions $d>1$, assuming rotational symmetry, we will consider the inverse of the \textit{radial} cdf of $f(t,\cdot)$ (i.e.\;of the partial mass function)---appropriately normalised. 
As in the first part of~\cite{carrillo_finite-time_2019}, we consider our equations posed on a bounded domain, more precisely on the centred ball $B(0,R_1)$ of radius $R_1>0$ with zero-flux boundary conditions.

\subsection{Change of variables}
  
 \subsubsection{One-dimensional case}\label{sssec:1d}
 Here, we consider the case $d=1$ and assume that $\gamma>2$, which represents the $L^1$-supercritical regime.
 The total mass of the initial datum $f_0$ is denoted by $m$.
 	Then, the equation satisfied by the pseudo-inverse $u(t,\cdot)$  of the \textit{cumulative distribution function} (cdf)
 	\begin{align*}
 	M(t,v)=\int_{-R_1}^v f(t,w)\,\d w,\quad v\in[-R_1,R_1],
 	\end{align*}
 	of $ f(t,\cdot)$  (cf.~\eqref{eq:defPsI}) formally states 
 	\begin{align}\label{eq:equ}
 	\partial_t u=(\partial_x u)^{-2}\partial_x^2 u-u(1+(u_x)^{-\gamma}), 
 	\end{align}
 	where $x\in(0,m)$ denotes the mass variable.
 	Formally, $f$ is related to $u$ by the identity
 	\begin{align*}
 	\partial_xu=\frac{1}{f(u)}.
 	\end{align*}
 	Upon multiplying eq.~\eqref{eq:equ} by the factor $(\partial_x u)^\gamma$, it can be rewritten as
\begin{align}\label{eq:invBefp1D}
 (\partial_x u)^{\gamma}\partial_t u-\frac{1}{\gamma-1}\partial_x\left((\partial_x u)^{\gamma-1}\right)+u((\partial_x u)^{\gamma}+1)=0.
\end{align}
While these new coordinates are generally known to be (numerically) favourable when investigating mass concentration phenomena in $1$D, a particular feature of equation~\eqref{eq:invBefp1D} is that the function $u\equiv0$, which at the level of $f$ corresponds to a Dirac delta at the origin, is an actual solution. Since mass conservation is a crucial feature of our Fokker--Planck model, the natural boundary condition for eq.~\eqref{eq:befp} states $\partial_rf+rf(1+f^\gamma)=0$ on $(0,\infty)\times\{-R_1,R_1\}$. It enforces the flux of particles through the boundary to be zero. Formally, at the level of $u$, this means that the RHS of eq.~\eqref{eq:equ} is zero on $(0,\infty)\times\{0,m\}$. Hence, if the solution $u$ is $C^{1,2}_{t,x}$ near and up to the boundary, this becomes $\partial_tu=0$ or, equivalently, 
\begin{align*}
  u=u_0\qquad\text{on }(0,\infty)\times\{0,m\}.
\end{align*}
This is the form we use in our numerical scheme. 
It corresponds to the Dirichlet conditions $u(t,0)=-R_1$, $u(t,m)=R_1$ for $t>0$.

As explained in Section~\ref{sec:introSim}, given a radius $R_1$ and a mass $m=\|f_0\|_{L^1(-R_1,R_1)}$ there exists a unique $\mu_\infty\in\mathcal{M}^+_b([-R_1,R_1])$ of mass $m$ which minimises the entropy $\mathcal{\widetilde H}$. At the level of $u$, we denote this minimiser by $u_\infty$. We further let $H(u):=\mathcal{H}(f)$ resp.~$\mathcal{\widetilde H}(\mu)$, where $\mu=u_\#\mathcal{L}^1$ is the push-forward measure  of the Lebesgue measure on $[0,m]$ under the map $u$ and will, in places, abbreviate $H_\infty:=H(u_\infty)=\mathcal{\widetilde H}(\mu_\infty)$. The dependence of $u_\infty$ on $R_1$ and $m$ will be omitted.
	For later reference, let us observe that $H(u)$ is  formally given by
\begin{align}\label{eq:Hu}
H(u)=\int_{(0,m)}  \left(\frac{|u|^2}{2}+\Psi(u_x)\right)\dd x,
\end{align}
where the function 
\begin{align}\label{eq:defPsi}
\Psi(s):=s\Phi(1/s)\quad\text{  is convex with }\quad \Psi''(s)=s^{-3}\Phi''(1/s)=\frac{1}{s^3 h(1/s)}.
\end{align}

 \subsubsection{Higher dimensions -- isotropic case}\label{sssec:eqhd}
 
 For isotropic solutions $f(t,v)=g(t,|v|),$ $v\in\mathbb{R}^d$, we can perform a similar transformation in higher dimensions. 
 In radial form, equation~\eqref{eq:invBefp1D} reads
\begin{align}\label{eq:radBefp}
  \partial_tg = r^{1-d}\partial_r\left(r^{d-1}\partial_rg+r^dg(1+g^\gamma)\right), \;t,r>0.
\end{align}
As a first ansatz one might try to consider the equation for the (pseudo-) inverse $R(t,z)$ of 
the \textit{radial cdf} $\bar M(t,r)=\int_0^{r}g(t,s)s^{d-1}\,\d s$. However, for bounded densities $f$ the function $\bar M$ is of class $O(r^d)$ as $r\to0$, implying that $R(t,\cdot)$ is at most $1/d$-H\"older near $z=0$ and $\partial_zR\gtrsim z^{1/d-1}\to\infty$ as $z\searrow0$, whenever $d>1$.
We therefore consider the normalised version $N(t,s)=\bar M(t,s^{1/d})$ or, equivalently,
\begin{align*}
  N(t,s)=\frac{1}{d}\int_0^{s}g(t,\sigma^{1/d})\,\d\sigma,
\end{align*}
which satisfies $\partial_sN(t,s)=\frac{1}{d}g(t,s^{1/d})$, and let $S(t,\cdot)$ denote the pseudo-inverse of $N(t,\cdot)$, so that $S=R^d$. From the formal relation $N(t,S(t,z))=z$ we deduce (omitting the time argument) 
\begin{align}\label{eq:f-S}
\partial_zS=\frac{d}{g(R)}.
\end{align}
Then, the equation~\eqref{eq:radBefp} for $g$ leads to the following equation for $S$:
\begin{align*}
  \frac{1}{d}\partial_tS - d\frac{S^{2-2/d}}{(\partial_zS)^2}\partial_z^2S+S(1+d^\gamma(\partial_zS)^{-\gamma})=0.
\end{align*}
Since we want our scheme to be able to deal with condensates, i.e.~$S(t,\cdot)\equiv0$ on some subinterval $(0,z(t))$, we multiply this equation by $(\partial_zS)^\gamma$ to obtain
\begin{align}\label{eq:Snonreg}
  (\partial_zS)^\gamma\tfrac{1}{d}\partial_tS - d\cdot S^{2-2/d}(\partial_zS)^{\gamma-2}\partial_z^2S+S((\partial_zS)^{\gamma}+d^\gamma)=0.
\end{align}
Notice that if $\gamma\in[1,2)$, the viscosity term has a factor which becomes unbounded when $S$ forms a condensate. We therefore consider 
for a small parameter $0<\varepsilon\ll1$
the following regularisation 
\begin{align*}
  (\partial_zS)^\gamma\tfrac{1}{d}\partial_tS - d\cdot S^{2-2/d}(\partial_zS+\varepsilon)^{\gamma-2}\partial_z^2S+S((\partial_zS)^{\gamma}+d^\gamma)=0
\end{align*}
or, equivalently, 
\begin{align*}
  \begin{cases}
    (\partial_zS)^\gamma\tfrac{1}{d}\partial_tS - \frac{d}{\gamma-1}\cdot S^{2-2/d}\frac{\d}{\d z}(\partial_zS+\varepsilon)^{\gamma-1}+S((\partial_zS)^{\gamma}+d^\gamma) =0, \;\;&\text{ if }\gamma>1,\\[3mm]
  (\partial_zS)^\gamma\frac{1}{d}\partial_tS - d\cdot S^{2-2/d}\frac{\d}{\d z}\log(\partial_zS+\varepsilon)+S((\partial_zS)^{\gamma}+d^\gamma) =0,\;\;&\text{ if }\gamma=1.
  \end{cases}
\end{align*}
We are mostly interested in the KQ model (where $\gamma=1$), and will thus focus on the equation
\begin{align*}
  d^{-1}\partial_zS\partial_tS - d S^{2-2/d}\frac{\d}{\d z}\log(\partial_zS+\varepsilon)+S(\partial_zS+d)=0,
\end{align*}
where $d=2,3$. Notice that a positive $\varepsilon$ decreases the strength of diffusion significantly when $\partial_zS\lesssim \varepsilon$. 
In order to counterbalance this effect, which may potentially lead to numerical artefacts when investigating the expected phenomenon of condensation, we  propose an artificial viscosity type regularisation of the form 
\begin{align}\label{eq:2regS3D}
  d^{-1}\partial_zS\partial_tS - d (S+\delta)^{2-2/d}\frac{\d}{\d z}\log(\partial_zS+\varepsilon)+S(\partial_zS+d)=0,
\end{align}
where $0<\delta\ll1$ is a small parameter.  
Below $\bar m$ (resp.\;$\bar m_c$) denotes the total mass of the initial datum $f_0$ (resp.\;of $f_c$) on $B(0,R_1)$ \textit{multiplied by the factor} $\frac{1}{|\partial B(0,1)|}$.
 Then, as in the $1$D case, the appropriate boundary conditions for equation~\eqref{eq:2regS3D} are $S(t,0)=0$ and $S(t,\bar m)=R_1^d$.

 As in Section~\ref{sssec:1d} we denote by $S_\infty=S_\infty(R_1,\bar m)$ the pseudo-inverse normalised radial cdf of the unique (isotropic) minimising measure in $\mathcal{M}^+_b(\overline{B}(0,R_1))$ corresponding to the choice $(R_1,m)$ of parameters, and generally let $H_d(S):=\mathcal{\widetilde H}(\mu)$, where $\mu$ is the unique isotropic measure in $\mathcal{M}^+_b(\overline{B}(0,R_1))$ satisfying $\mu(\overline{B}(0,r))=\nu([0,r^d])\cdot |\partial B(0,1)|$ and $\nu$ denotes the measure associated with the generalised inverse of $S$. We also abbreviate $H_\infty:=H(S_\infty)$ and $H(t):=H_d(S(t))$.

\subsubsection{Higher dimensions -- anisotropic case}	

Let us briefly discuss that one can perform a related change of variables in higher dimensions without radial symmetry. With this aim, one needs to consider vector-valued transformations $u(t,\cdot):U\to V$,  $U,V\subset\mathbb{R}^d$, 
which are formally related to the original density $f$ via 
\begin{align*}
\det \nabla u(t,x)\cdot f(t,u)=1.
\end{align*}
Similarly to~\cite{evans_diffeo_2005,carrillo_moll_2009}, one finds that the system governing the evolution of $u=(u^1, \dots, u^d)^T$ can formally be written as
\begin{align*}
\quad\left[(\det \nabla u)^2\Psi''(\det \nabla u)\right]\partial_tu^i -\partial_{x_k}\left( \Psi'(\det \nabla u)(\mathrm{cof}(\nabla u))_k^i\right)  + u^i = 0
\end{align*}
for $i=1,\dots,d$, where $\Psi$ is defined as in~\eqref{eq:defPsi}. 
The entropy $H_{\mathrm{ani},d}(u)$ in the new variables takes the form
\begin{align*}
H_{\mathrm{ani},d}(u) = \int_U\left(\tfrac{1}{2}|u|^2+\Psi(\det\nabla u)\right)\,\dd x.
\end{align*}
Observe that in the vectorial case $H_{\mathrm{ani},d}(u)$ is no longer convex but merely polyconvex in $\nabla u$. 
This route could potentially allow to numerically analyse concentrations without radial symmetry in higher dimensions, 
as it is the case in 2D for aggregation and Keller--Segel type problems close to the blow-up time~\cite{carrillo_numerical_2016}.
Even if this method deserves further exploration, we focus here on the isotropic case to capture the direct generalisation of the 1D behaviour in the 3D realistic setting.

\subsection{The semidiscrete scheme}\label{ssec:semidiscrete}

The scalar equations~\eqref{eq:invBefp1D} and~\eqref{eq:2regS3D} are discretised fully implicitly in time.
We let $\tau$ be the discrete time step and denote by $\{u^{n}\}_{n\in\mathbb{N}}$ the time-discrete solution of the implicit Euler discretisation of equation~\eqref{eq:invBefp1D}. More precisely, given a non-decreasing function $u^{n}$ satisfying $u^{n}(0)=-R_1$ and $u^{n}(m)=R_1$,
 the problem for $u=u^{n+1}$ reads
\begin{align}\label{eq:invBefp1Dsemi}
  \left(\partial_x u\right)^\gamma \frac{u-u^{n}}{\tau} 
  -\tfrac{1}{\gamma-1}\partial_x\left((\partial_x u)^{\gamma-1}\right)+u((\partial_x u)^{\gamma}+1)=0
\end{align}
subject to the Dirichlet boundary conditions $u^{n+1}(0)=-R_1, u^{n+1}(m)=R_1$.

Let us here make a short digression to explain the main difference and potential difficulty of the present problem with respect to the Wasserstein  gradient flows treated in \cite{blanchet_convergence_2008,carrillo_numerical_2016}.
Those works are based on the idea that the Wasserstein gradient flow  of the entropy/free energy in the original variables is equivalent to an $L^2$ gradient flow for the problem in the $u$-variables.
Loosely speaking, the semidiscrete $L^2$ gradient flow for $H(u)$ reads as follows: given $\tilde u^n$ formally define $\tilde u^{n+1}$ as a solution of the problem 
\begin{align*}
\tilde u^{n+1}\in\mathrm{arg}\inf_{\tilde u}\left\{\frac{1}{2\tau}\|\tilde u-\tilde u^n\|_{L^2}^2+H(\tilde u)\right\}.
\end{align*}
The associated  Euler--Lagrange equations, $\frac{\tilde u-\tilde u^n}{\tau}=-\nabla_{L^ 2} H(\tilde u),$ read
\begin{align*}
\frac{\tilde u-\tilde u^n}{\tau}=-[-\partial_x(\Psi'(\tilde u_x))+\tilde u].
\end{align*}
To compare this with our problem, 
we write eq.~\eqref{eq:invBefp1Dsemi} in the more concise equivalent form
\begin{align*}
u_x^2\Psi''(u_x) \frac{u-u^n}{\tau}=-[-\partial_x(\Psi'(u_x))+u],
\end{align*}
which suggests that in some sense a gradient flow structure is kept. 
At least, as will be shown below, we keep an important property in the semidiscrete numerical scheme, namely the monotonicity of the entropy. 
Recall that in 1D the entropy  $H(u)$ in the $u$-variables (see~\eqref{eq:Hu}) is convex in the classical sense, and
 it is well-known that the implicit Euler scheme applied to a gradient flow of a convex functional satisfies the semidiscrete entropy inequality $H(\tilde u^{n+1})\leq H(\tilde u^{n})$ for all $n$.  
In our situation, thanks to the convexity of the integrand of $H$, the entropy decay along the sequence $\{u^n\}$ can be recovered by a simple estimate:
\begin{align*}
H(u)-H(u^n)&\le \int_{(0,m)} (u(u-u^n)+\Psi'(u_x)(u-u^n)_x)\,\dd x
\\&= \int_{(0,m)}  (u-\partial_x(\Psi'(u_x)))(u-u^n)\,\dd x
\\&=-\tau\int_{(0,m)}  u_x^2\Psi''(u_x)\Big|\frac{u-u^n}{\tau}\Big|^2\,\dd x\le0.
\end{align*}
Here, we used the fact that in the above integration by parts the boundary terms vanish since, by construction,  $u=u^{n}$ on $\partial(0,m)$. This shows the entropy decay property of the semidiscrete scheme~\eqref{eq:invBefp1Dsemi}: $H(u^{n+1})\leq H(u^{n})$ for all $n$. We note that similar properties with a similar strategy of proof are found for finite volume schemes of gradient flows \cite{BCH}.

\begin{remark}[Higher dimensions, isotropic case]
In higher dimensions the entropy $H_d(S)$, introduced in Section~\ref{sssec:eqhd}, takes the form (see also~\eqref{eq:f-S})
\begin{align}\label{eq:HdS}
H_d(S) = \int \left(\tfrac{1}{2}S^\frac{2}{d}+ \Psi_d(\partial_zS)\right)\,\dd z,
\end{align}
where $ \Psi_d(s) =\Psi(\tfrac{s}{d})$ is again convex. 
If $d=2$, thanks to convexity, the implicit Euler discretisation of eq.~\eqref{eq:Snonreg} can be shown to keep the entropy decay by arguing as in the 1D case. In higher dimensions, $d>2$, this argument breaks down due to the kinetic part of the entropy failing to be a convex function of $S$. Notice, however, that the convexity in the highest order term, $\partial_zS$, is maintained.
\end{remark}

\subsection{The fully discrete scheme}

The semidiscrete nonlinear system \eqref{eq:invBefp1Dsemi} is discretised using finite differences and solved by the Newton--Raphson method. 
In the one dimensional case, the finite difference approximation in space is chosen in such a way as to preserve the equation's symmetry, viz. 
\small
\begin{align}\label{eq:discreteu}
  (u_{i+1}^n-u_{i-1}^n)^\gamma(2h)^{-\gamma} \frac{u_i^n-u^{n-1}_i}{\tau} 
  - ((u_{i+1}^n-u_{i}^n)^{\gamma-1}-(u_{i}^n-u_{i-1}^n)^{\gamma-1})h^{-\gamma}(\gamma-1)^{-1}&
  \\\nonumber +u_i^n((u_{i+1}^n-u_{i-1}^n)^\gamma(2h)^{-\gamma}+1)=0,&
\end{align}
\normalsize
for $i=1,\dots,N-1,$ complemented with the boundary conditions $u^n_0=u^0_0=-R_1$ and $u^n_N=u^0_N=R_1$.
We use a similar full discretisation for~\eqref{eq:2regS3D}, viz.
\begin{multline*}
  (S^n_{i+1}-S^n_{i-1})(2hd\tau)^{-1}(S^n_i-S^{n-1}_i)\\- d (S_i^n+\delta)^{2-2/d}(\log((S^n_{i+1}-S^n_{i})/h+\varepsilon)-\log((S^n_{i}-S^n_{i-1})/h+\varepsilon))/h
  \\+S^n_i((S^n_{i+1}-S^n_{i-1})/(2h)+d)=0
\end{multline*}
for $i=1,\dots,N$, 
where the boundary conditions are given by $S_0^n=S_0^0=0$ and $S_N^n=S_N^0=R_1^d$.

\paragraph{Algorithm.}

Given $u^{n-1}$ the discrete approximation $u^{n}$ at the subsequent time point is computed using a Newton--Raphson iteration. The iteration is stopped as soon as the smallness condition $\|F_{\mathrm{NR}}(u^n,u^{n-1},h,\tau)\|_{l^2}<10^{-8}$ is satisfied, where $F_{\mathrm{NR}}(u^n,u^{n-1},h,\tau)_i$ is given by the LHS of equation~\eqref{eq:discreteu} multiplied by $h^\gamma$. For $S$ we proceed similarly.

\begin{remark}
In the simulations exhibiting the numerically somewhat delicate condensation phenomenon, the inverse cdf becomes slightly non-monotonic during the Newton--Raphson iteration,  which leads to very small imaginary parts in the above scheme and of the solution at the subsequent time step. In our actual code we therefore rearrange the approximation in each Newton--Raphson iteration to ensure monotonicity. Alternatively, one can replace the first derivatives $u_x$ by their absolute values $|u_x|$ and discretise and simulate this equation.
In practice, the differences between the results using the first and the second option are negligible.
A similar statement applies to the higher-dimensional case, where we choose again the option of the monotonic rearrangement.
\end{remark}

 \section{Numerical experiments}
 In this section we describe the validation of our scheme, and present and discuss our numerical experiments.
 
  \subsection{\texorpdfstring{$L^1$}{}-supercritical bosonic Fokker--Planck model in  1D:    \texorpdfstring{\\}{}
 simulations replicating the theory}
 \label{ssec:sim1D}
 
First, we demonstrate the reliability of the proposed numerical scheme in~1D by reproducing the features proved in~\cite{carrillo_finite-time_2019}. 
In addition,  we use the scheme to predict that the entropy decays at an exponential rate, even after the onset of a condensate.

 If not stated otherwise, we choose $\gamma=2.9$ and use a centred Gaussian as initial datum, viz.
 \begin{align}\label{eq:init}
   f_0(v) = A\mathrm{e}^{-\frac{|v|^2}{2\sigma^2}}
 \end{align}
 for fixed positive constants $A$ and $\sigma$. Moreover, we always set $R_1=1$.
 We remark that for $d=1$ and the above choice of $\gamma$ and $R_1$ the critical mass $m_c$ takes the numerical value $m_c\approx5.37$.

 \subsubsection{Validation in 1D} 
 
 We begin with validating the $1$D scheme~\eqref{eq:discreteu} by comparing the solution for a given mesh with a numerical reference solution calculated on a fixed and much finer mesh. 
 We set $\sigma=0.7$, $A=4.5$ in~\eqref{eq:init} as well as $T=0.025$. For simplicity, the mass variable $x\in[0,m]$ is often referred to as the \textit{spatial} variable. The numerical reference solution is computed on a grid of $12801$ (equidistant) spatial mesh points and a total number of $1000$ (equidistant) time points. Notice that the values of the parameters $A$ and $\sigma$ coincide with those in~\ref{it:sup} below and observe that, in the simulations based on~\ref{it:sup}, well before the final time $T=0.025$ chosen for our validation, a significant amount of mass has accumulated at the origin (cf.\;Figures~\ref{fig:111} and~\ref{fig:113}). Therefore, our validation covers the case in which condensation occurs.
 
\begin{table}[H]\footnotesize
\parbox{.49\linewidth}{
\centering
\begin{tabular}{|c c c c|}                
\hline 
 timesteps & meshsize & $L^2_{x}$ error& rate \\  
\hline  \hline                                  
1000 & 50 &  7.3825e-3 & - \\       
\hline                                    
1000 & 100 & 2.1290e-3  & 1.7939 \\ 
\hline                                    
1000 & 200 & 5.6056e-4 & 1.9253 \\ 
\hline                                    
1000 & 400 & 1.4222e-4 & 1.9788 \\ 
\hline                                    
1000 & 800 & 3.5598e-5 & 1.9982 \\ 
\hline                                    
1000 & 1600 &  8.8061e-6 & 2.0152 \\
\hline                                    
1000 & 3200 & 2.0991e-6 & 2.0687 \\
\hline                                
\end{tabular}   
     \caption[Validation w.r.t.\;reference at time $T=0.025$ ($d=1,\gamma=2.9,m>m_c$)]{Convergence to reference solution at time $T=0.025$.}  
      \label{table:1Dmesh} 
}
 \hfill
\parbox{.49\linewidth}{
\centering
      \begin{tabular}{|c c c c|} 
    \hline     
     timesteps & meshsize & $L^2_{t,x}$ error& rate \\  
\hline  \hline  
      10& 50&  6.1372e-3  & - \\  
      \hline  
      20& 100& 3.1393e-3  &  0.9671 \\ 
      \hline  
      40& 200& 1.5817e-3   &   0.9890 \\  
      \hline  
      80& 400& 7.8542e-4    & 1.0099 \\  
      \hline  
      160& 800& 3.8200e-4   &  1.0399 \\  
      \hline  
      320& 1600& 1.7877e-4  &   1.0955 \\  
      \hline  
      640& 3200& 7.6728e-5  &  1.2203 \\   
      \hline 
      \end{tabular}  
      \caption[Validation w.r.t.\;reference on space-time grid ($d=1,\gamma=2.9,m>m_c$)]{Convergence to reference solution (on space-time grid).}  
      \label{table:1Dfull} 
}
\end{table}
    
     Table~\ref{table:1Dmesh} displays the discrete $L^2_x$ error of the solution on the coarser mesh with respect to the reference solution, evaluated at the final time $T$, while
 Table~\ref{table:1Dfull} indicates the $L^2$ space-time error between computed and reference solution.
 The results suggest a second order dependence of the error on the spatial increment and a first order dependence on the temporal increment. 
 As long as the solution is not degenerate, this can be explained by the fact that
 we use an implicit Euler scheme in time (which is first-order accurate), a central finite difference discretisation in space (whose truncation error is of second order) and have chosen a high resolution in time for the test using purely spatial refinement, which makes the temporal error negligible in this test. 
 Notice, however, that the degenerate case requires more care and that, in this work, we do not provide a rigorous numerical analysis of the scheme.
 
 \begin{remark}
 Higher-order implicit time discretisations could be considered to obtain better accuracy.  Table~\ref{table:CN} displays the convergence rates upon refinement of the space-time mesh using a Crank--Nicolson-type (CN) time discretisation for \eqref{eq:invBefp1D} with parameters \ref{it:sub} and  clearly confirms the second order accuracy of CN. 
However, we would like to point out that the initial datum determined by~\ref{it:sub} is mass-subcritical, and therefore the 2nd order accuracy is obtained for smooth solutions. Our simulations beyond blow-up indicate that the Newton solver for the implicit Euler scheme has better stability properties to cope with condensates than the CN scheme though.
\end{remark}
\begin{table}[H]\footnotesize\centering
	\parbox{.49\linewidth}{
		\centering
		\begin{tabular}{|c c c c|} 
			\hline     
			timesteps & meshsize & $L^2_{t,x}$ error& rate \\  
			\hline  \hline  
			10& 50&  5.2392e-3  & - \\  
			\hline  
			20& 100& 1.1085 e-3  & 2.2408 \\ 
			\hline  
			40& 200& 2.4257 e-4   &   2.1921\\  
			\hline  
			80& 400& 5.6873e-05    & 2.0926 \\  
			\hline  
			160& 800& 1.3983e-05       &  2.0241 \\  
			\hline  
		\end{tabular}   
		\caption{Convergence to reference solutions using CN and \ref{it:sub}.}  
		\label{table:CN} 
	}
\end{table}
    
\subsubsection{Comparing simulations and theoretical results}

In order to numerically confirm the dynamical properties of eq.~\eqref{eq:befp} in $1$D, we run our scheme with the following four sets of parameters covering  the mass-super resp.\;-subcritical, the asymmetric case as well as the case of the initial datum being highly concentrated near the origin $v=0$:
\begin{enumerate}[label=(P\arabic*)]
 \item\label{it:sup} $m>m_c:$ $\sigma=0.7$, $A=4.5$, $T=0.4$, $\tau=0.001$, $n=2001$ ($n:=$ number of spatial grid points).
 \item\label{it:supasy} Asymmetric \& $m>m_c:$ translated Gaussian 
 $f_0(v) = A\mathrm{e}^{-|v-v_0|^2/(2\sigma^2)}+0.1$ chosen as initial datum using the parameters $v_0=-1$,  $\sigma=0.7$ and $A=4.5$. Moreover, $T=0.4$, $\tau=0.001$, $n=2001$. The shift by $+0.1$ ensures that the cdf of $f_0$ is numerically still well invertible close to $v=R_1$.
 \item\label{it:sub} $m<m_c:$  $\sigma=0.7$, $A=1.5$, $T=0.4$, $\tau=0.001$, $n=2001$.
 \item\label{it:subbu} Concentrated \& $m<m_c:$  $\sigma=0.1$, $A=1.5$, $T=0.4$, $\tau=10^{-6}$, $n=10001$.
\end{enumerate}

The approximate total mass for each of these simulations is indicated in part~\textbf{(a)} of the corresponding figure: it is the maximal value of the part of the horizontal axis which is displayed.

\begin{figure}[H]
\begin{minipage}{0.49\textwidth}\centering
\includegraphics[scale=\sca]{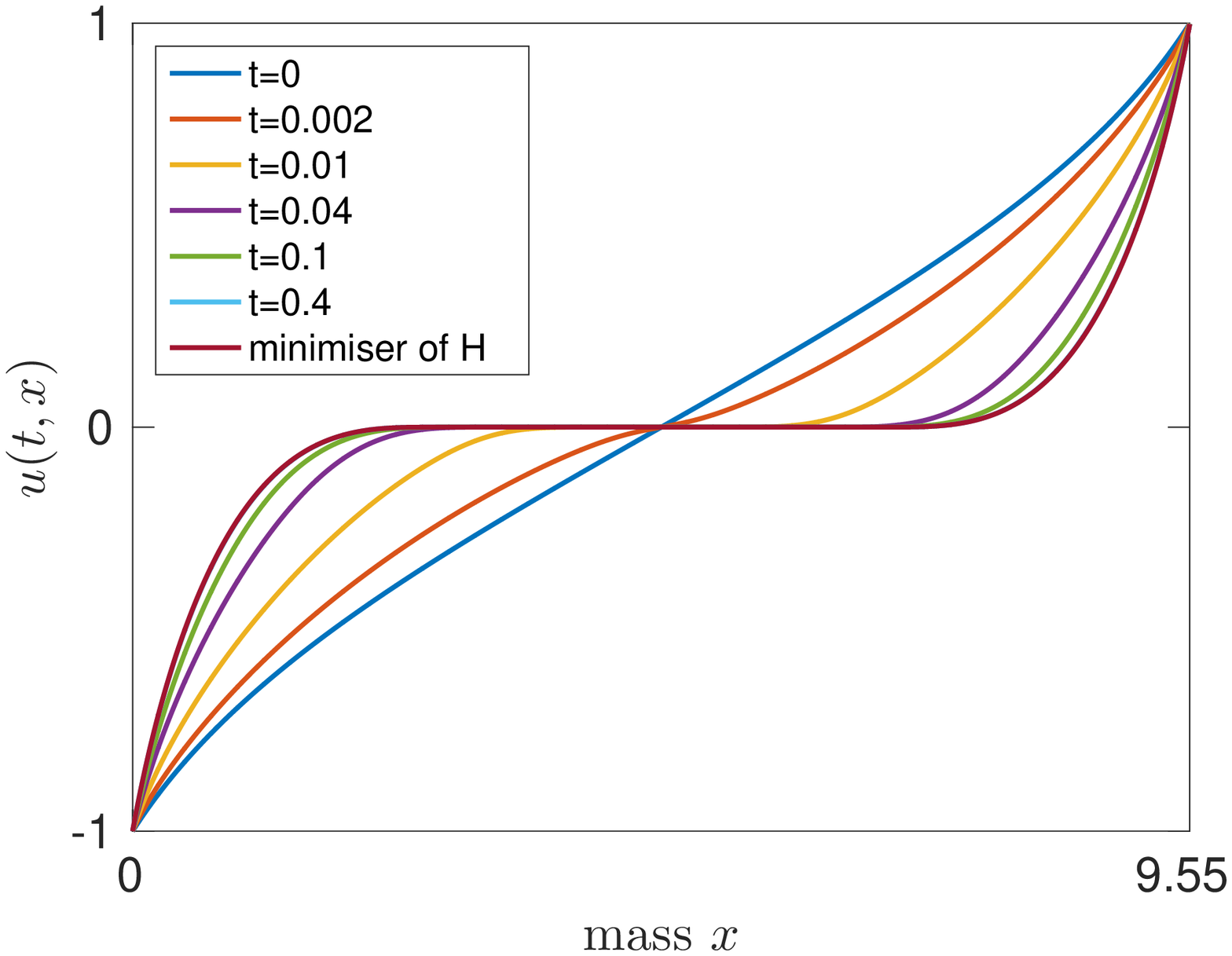} 
\subcaption{$u(t,\cdot)$ and $u_\infty$.}
\label{fig:111}
\end{minipage} 
\begin{minipage}{0.49\textwidth}\centering
\includegraphics[scale=\sca]{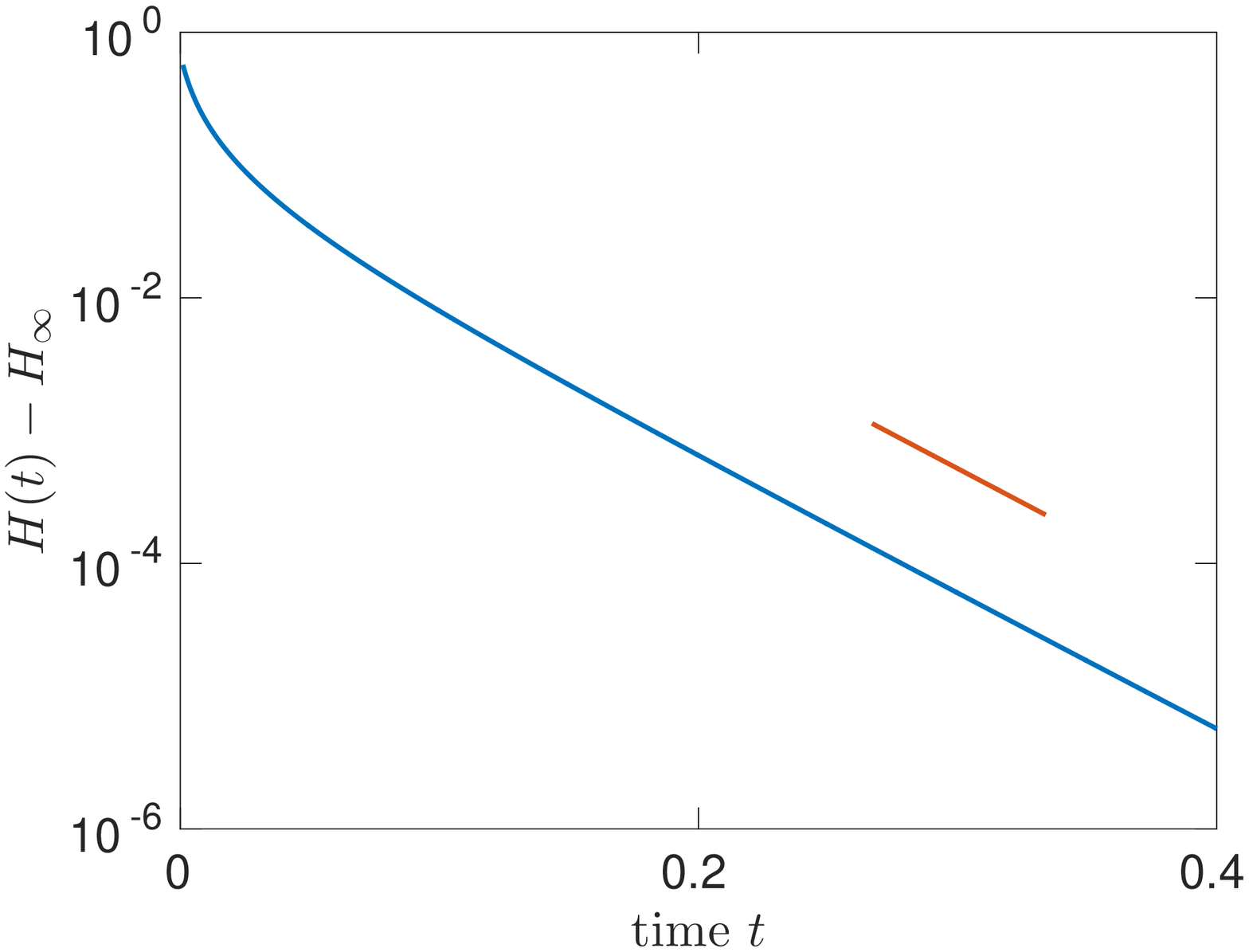}
\subcaption{Evolution of the relative entropy.}
\label{fig:112}
\end{minipage}

\vspace{.5cm}

\begin{minipage}{0.49\textwidth}\centering
\includegraphics[scale=\sca]{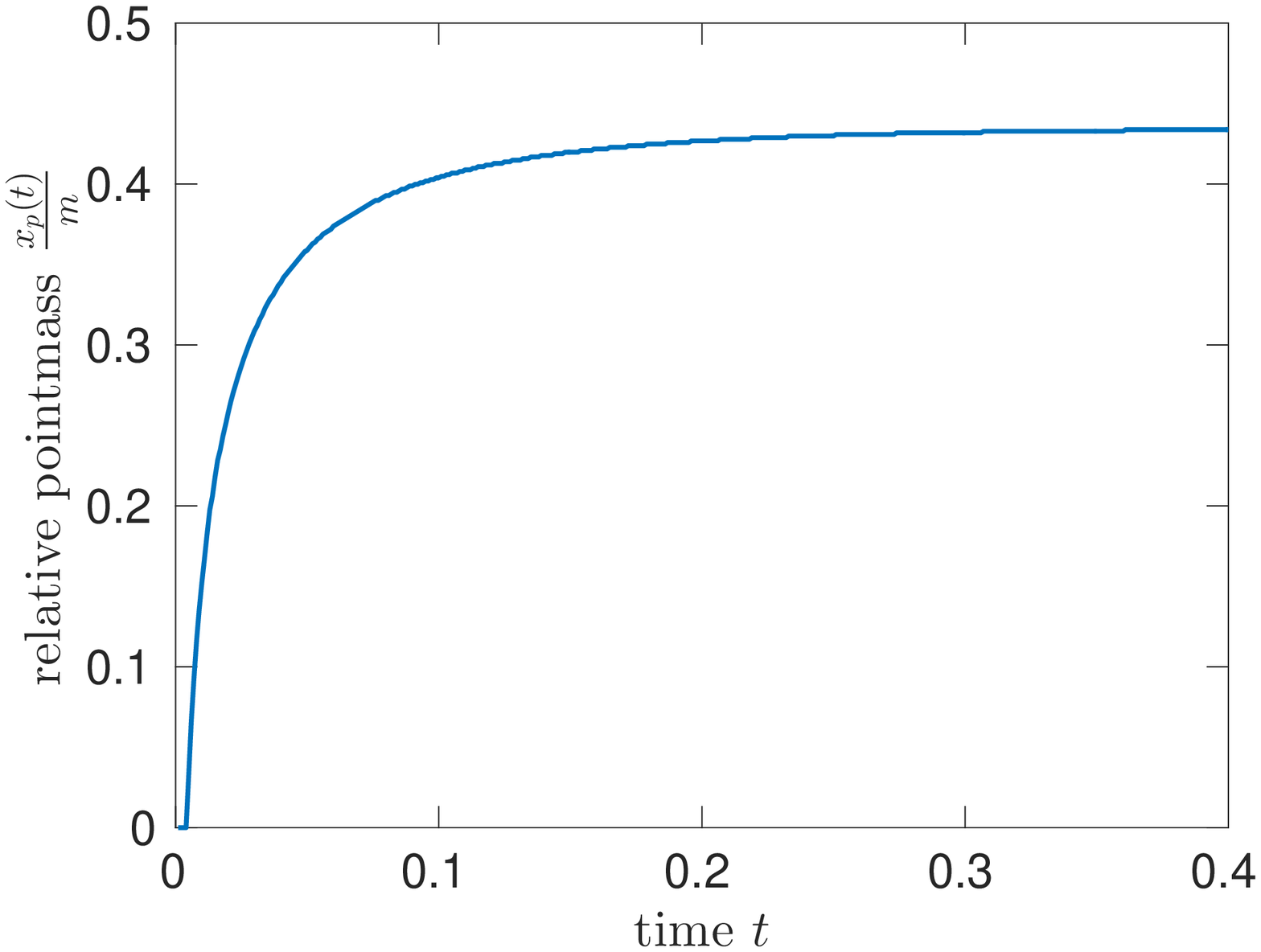} 
\subcaption{Evolution of the Dirac part.}
\label{fig:113}
\end{minipage}
\begin{minipage}{0.49\textwidth}\centering
\includegraphics[scale=\sca]{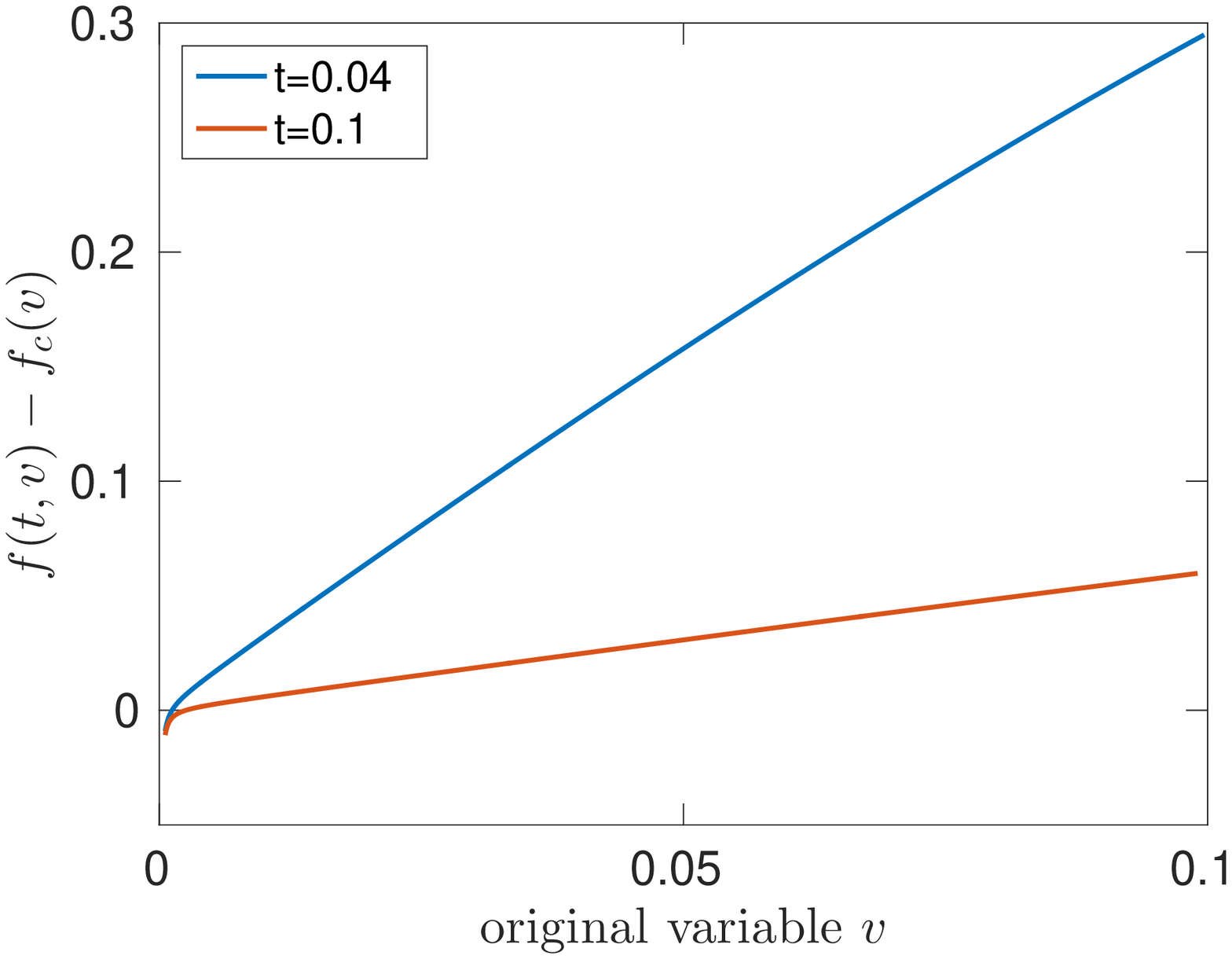} 
\subcaption{Behaviour near singularity.}
\label{fig:115}
\end{minipage}
\caption[Large-time behaviour for~\ref{it:sup} ($d=1, \gamma=2.9, m>m_c$)]{Long-time behaviour in the mass-supercritical case~\ref{it:sup}  ($d=1, \gamma=2.9$).}
\end{figure}

\begin{figure}
\begin{minipage}[H]{0.49\textwidth}\centering
\includegraphics[scale=\sca]{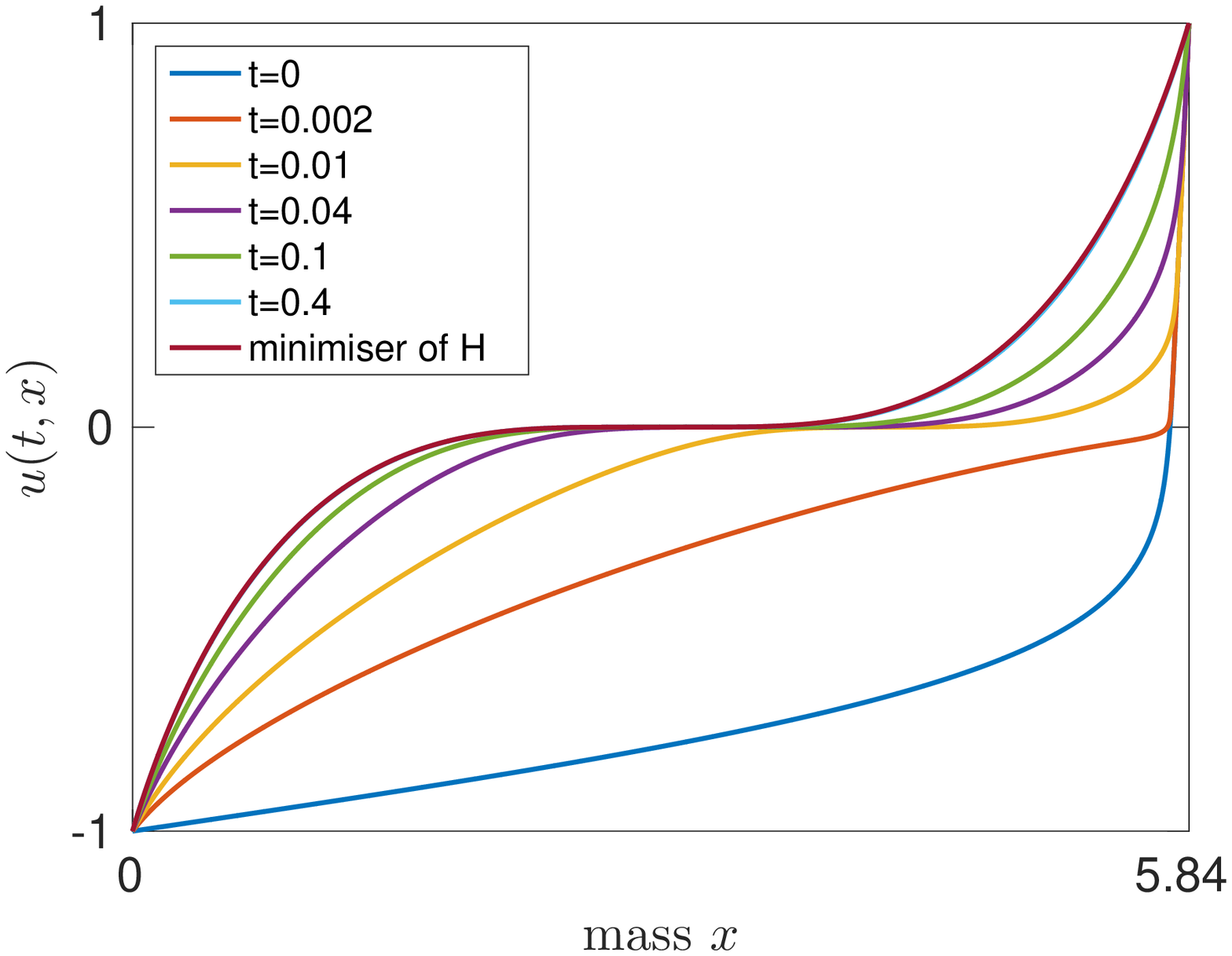} 
\subcaption{$u(t,\cdot)$ and $u_\infty$.}
\label{fig:141}
\end{minipage} \hfill 
\begin{minipage}[H]{0.49\textwidth}\centering
\includegraphics[scale=\sca]{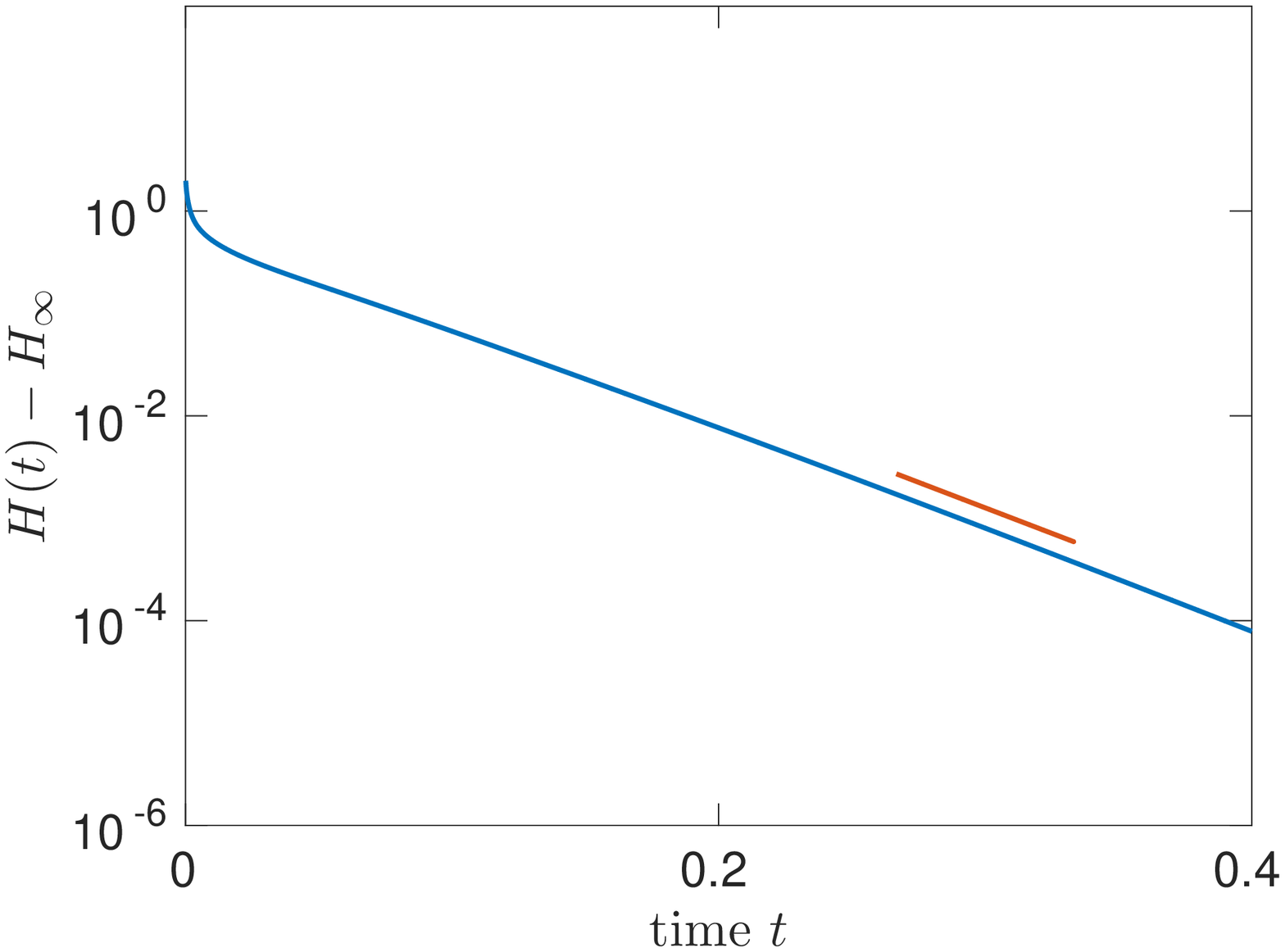}
\subcaption{Evolution of the relative entropy.}
\label{fig:142} 
\end{minipage}

\vspace{.5cm}

\begin{minipage}[H]{0.49\textwidth}\centering
\includegraphics[scale=\sca]{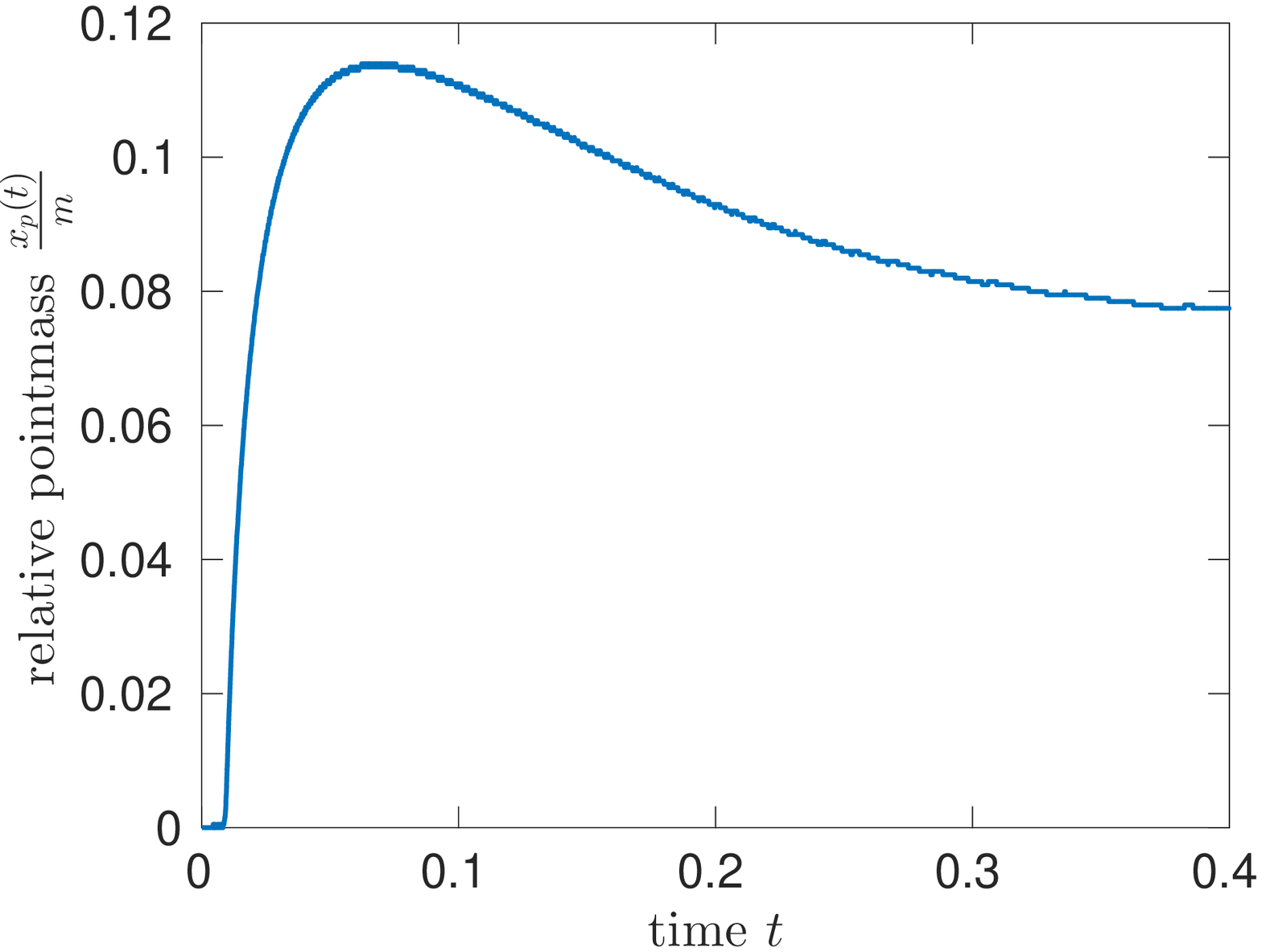} 
\subcaption{Evolution of the Dirac part.}
\label{fig:143}
\end{minipage} \hfill 
\begin{minipage}[H]{0.49\textwidth}\centering
\includegraphics[scale=\sca]{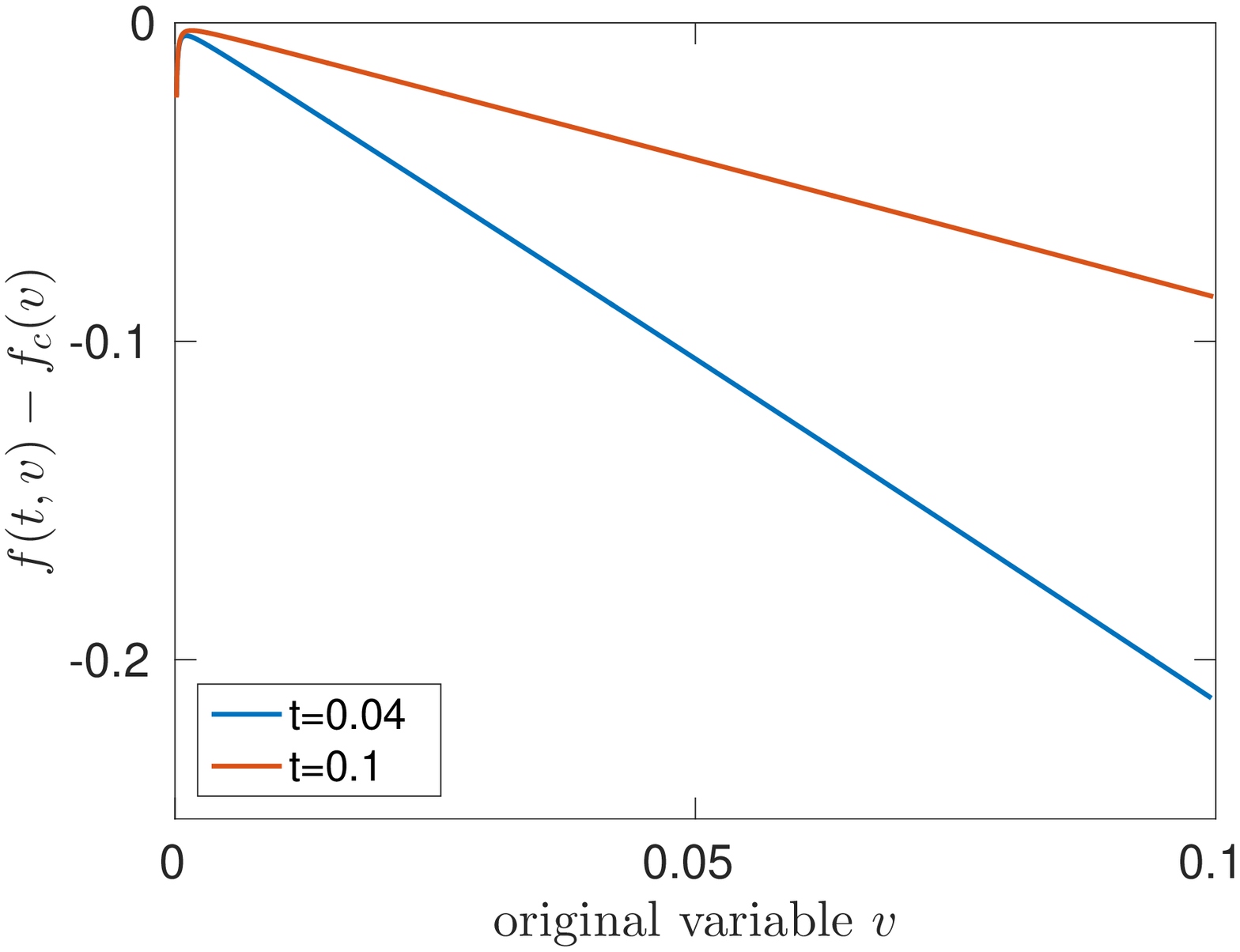} 
\subcaption{Behaviour near singularity.} 
\label{fig:145}
\end{minipage}
\caption[Long-time behaviour for~\ref{it:supasy} ($d=1, \gamma=2.9, m>m_c$, asymmetric)]{Long-time behaviour for asymmetric mass-supercritical datum~\ref{it:supasy} ($d=1, \gamma=2.9$).
}
\label{fig:image2}
\end{figure}

\begin{figure}[ht!]
 \begin{minipage}{.49\textwidth}\centering
   \includegraphics[scale=\sca]{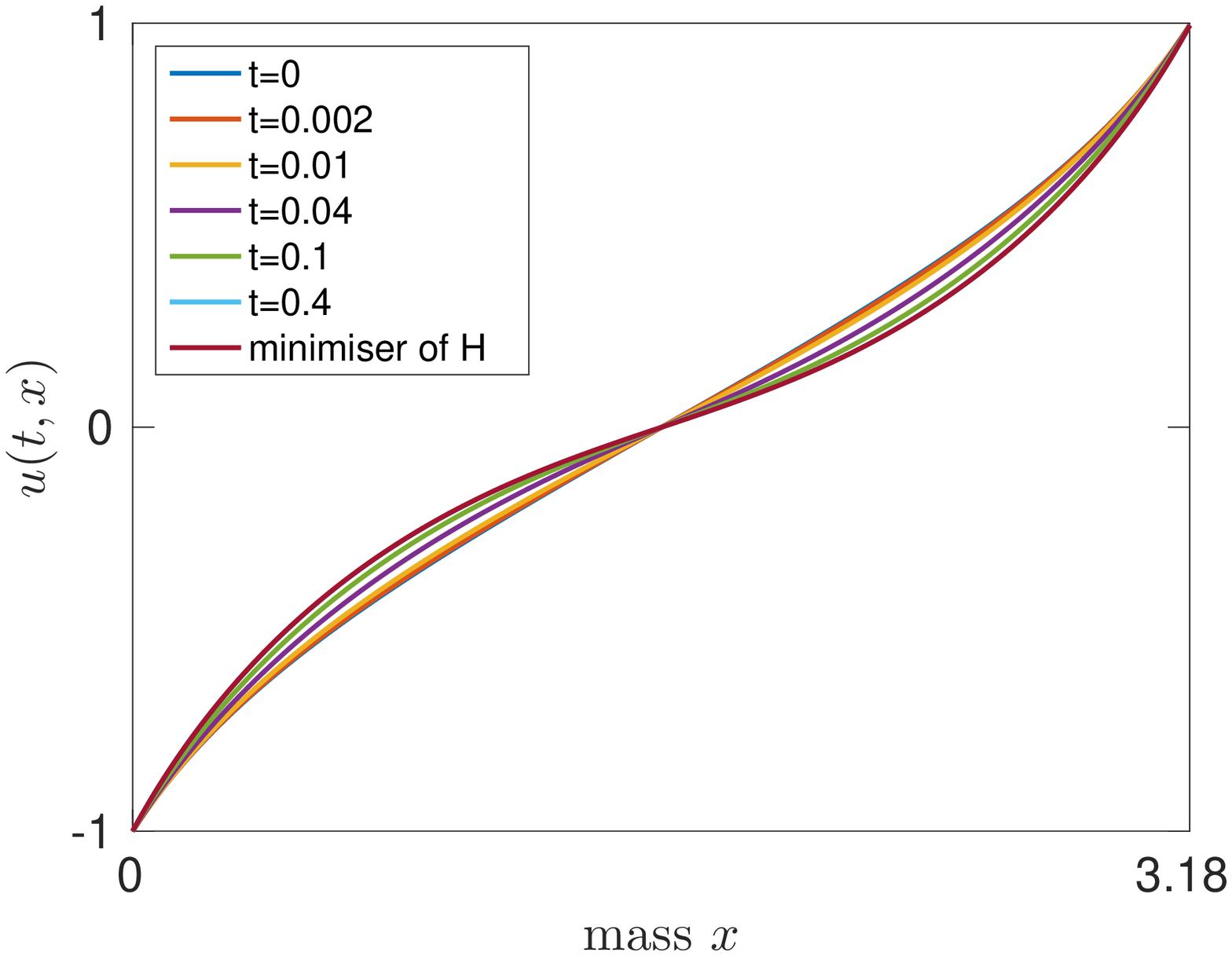}
     \subcaption{$u(t,\cdot)$ and $u_\infty$ ($\sigma=0.5$).}
\label{fig:121}
\end{minipage}
\begin{minipage}{.49\textwidth}\centering
\includegraphics[scale=\sca]{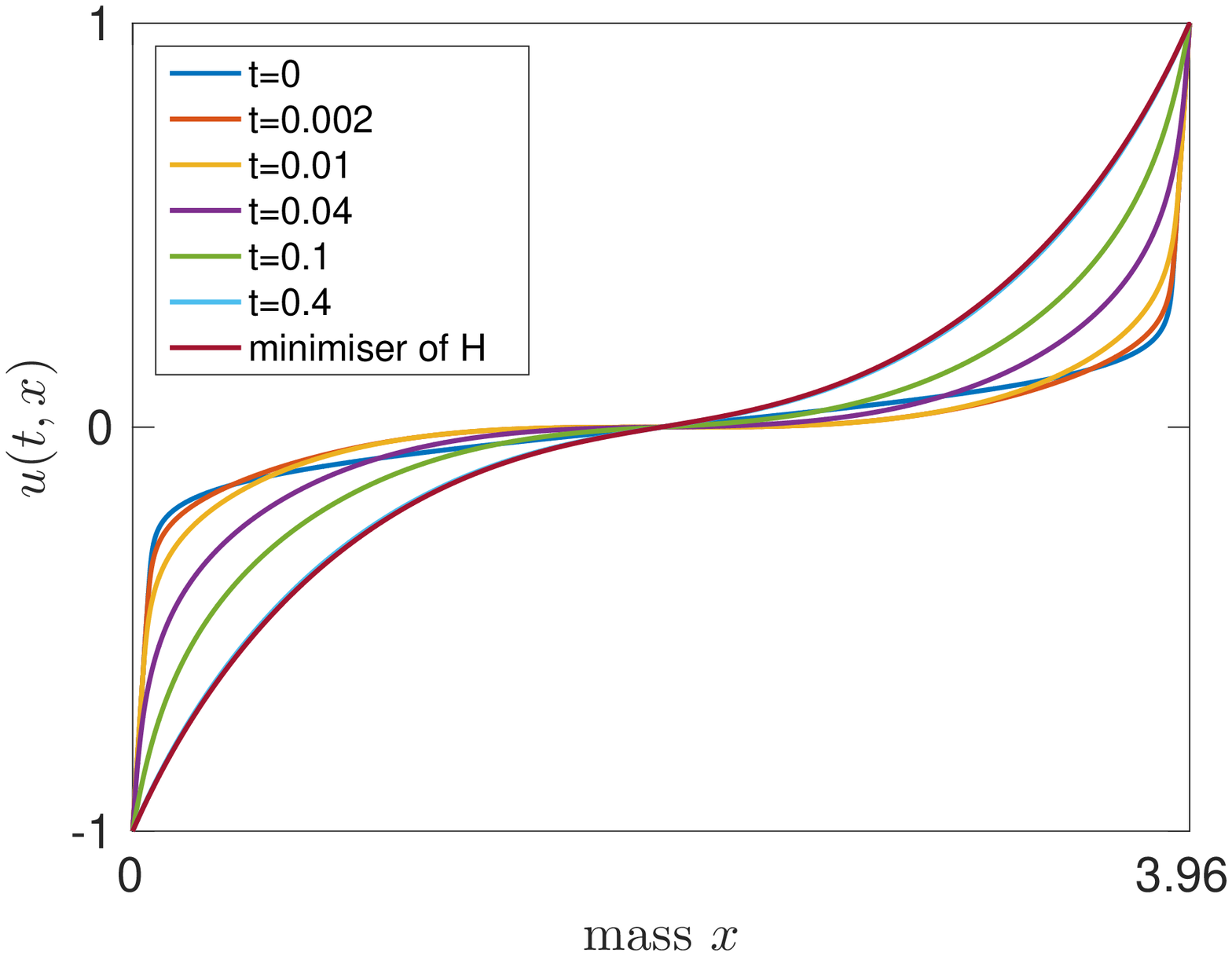} 
\subcaption{$u(t,\cdot)$ and $u_\infty$ ($\sigma=0.1$).}
\label{fig:131}
\end{minipage}
\vspace{0.5cm}

\begin{minipage}{.49\textwidth}\centering
\includegraphics[scale=\sca]{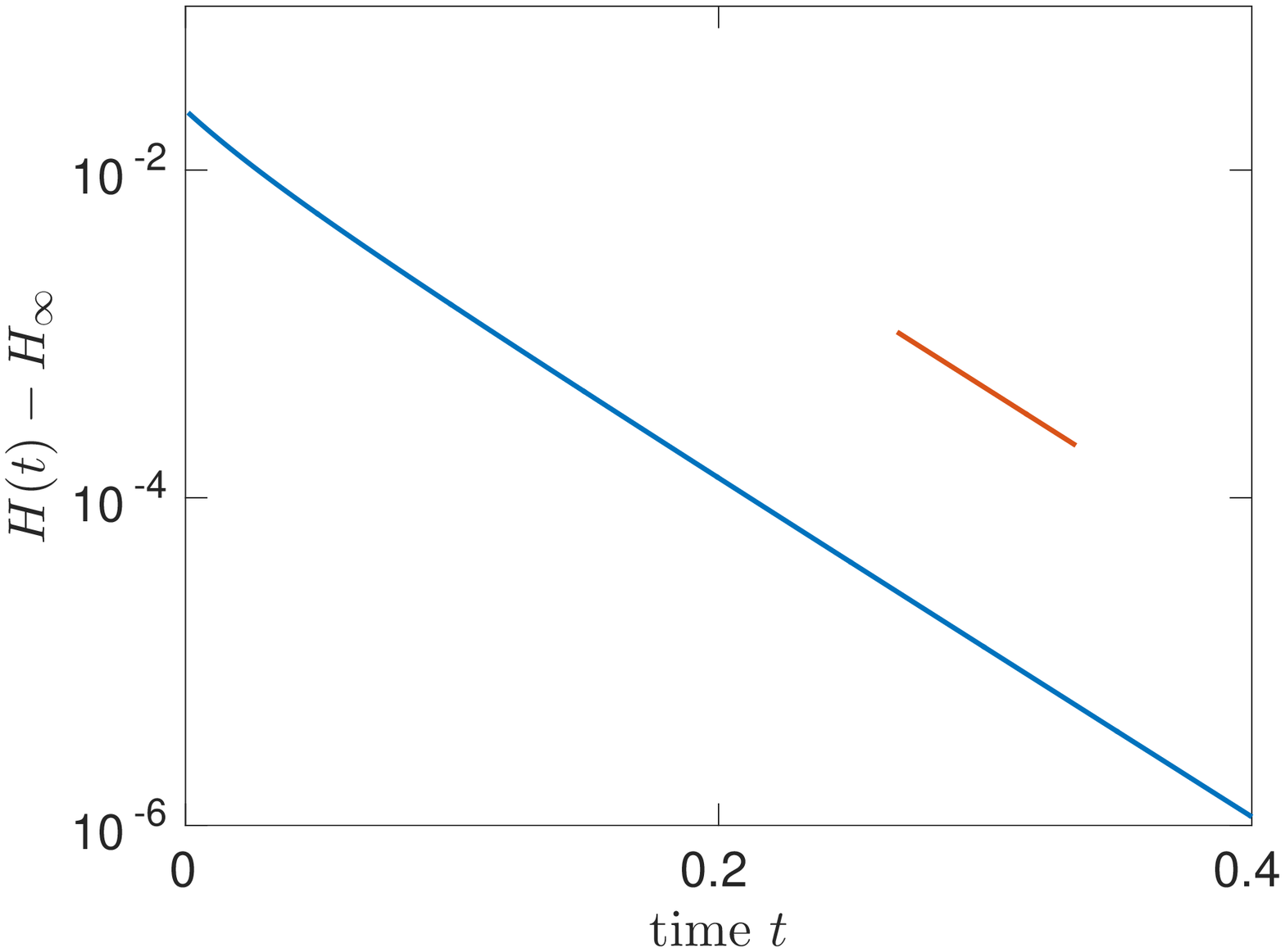}
\subcaption{Evolution of the relative entropy \\($\sigma=0.5$).}
\label{fig:122}
\end{minipage}
\begin{minipage}{.49\textwidth}\centering
\includegraphics[scale=\sca]{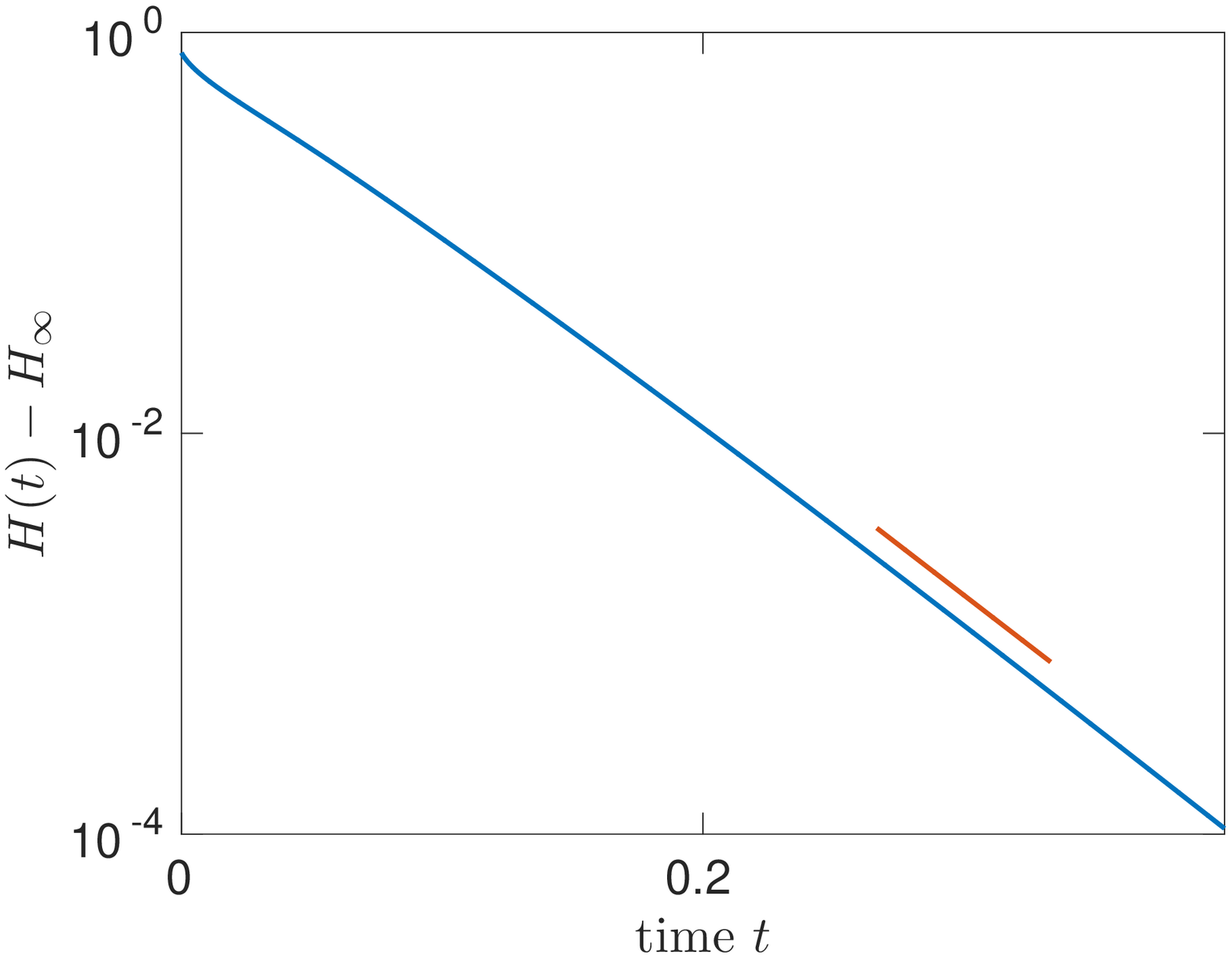} 
\subcaption{Evolution of the relative entropy \\($\sigma=0.1$).}
\label{fig:132}
\end{minipage}
\vspace{0.5cm}

\begin{minipage}{.49\textwidth}\centering
\includegraphics[scale=\sca]{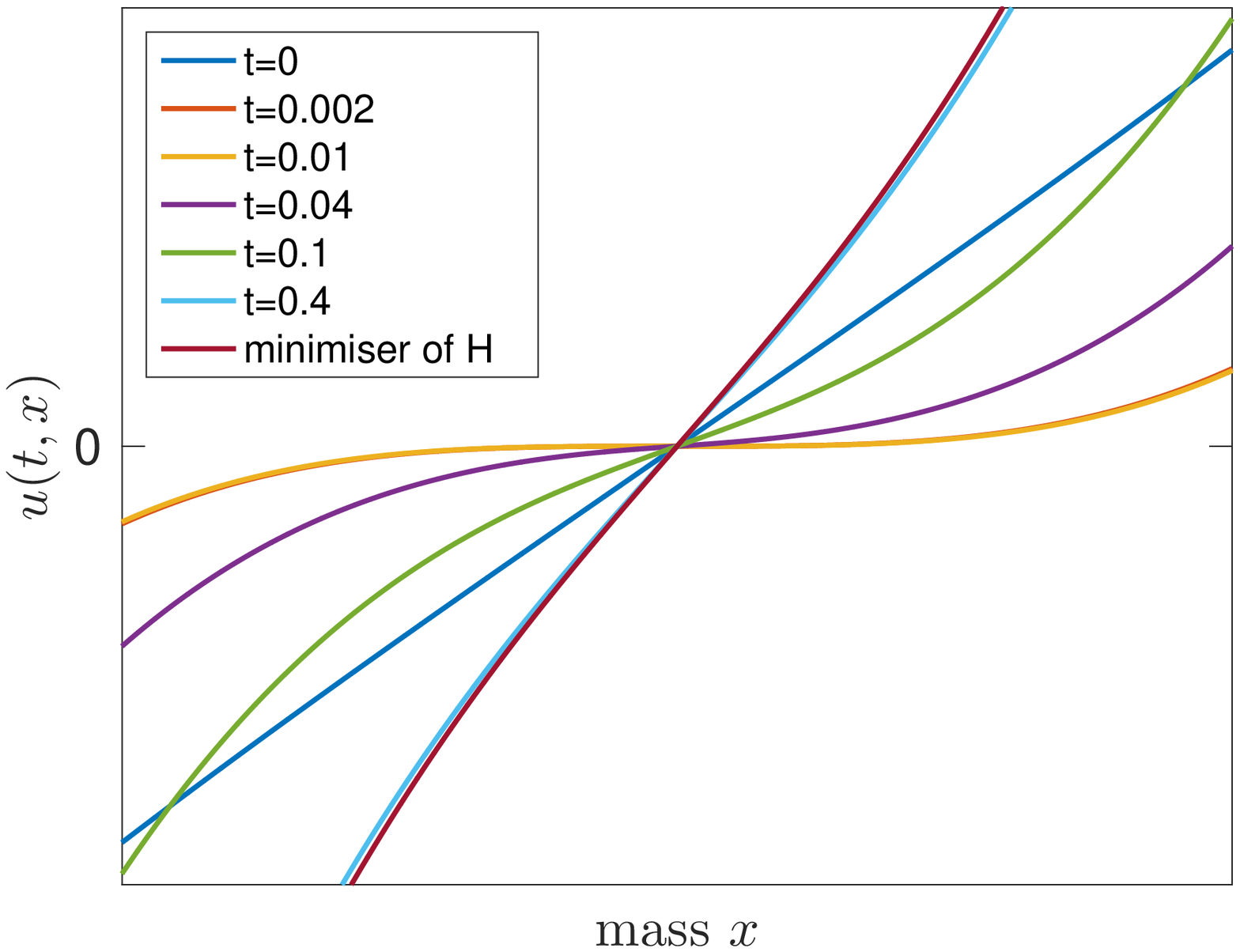} 
\subcaption{
Zoomed-in view of Fig.~\ref{fig:131}.}
\label{fig:135}
\end{minipage}
\begin{minipage}{.49\textwidth}\centering
\includegraphics[scale=\sca]{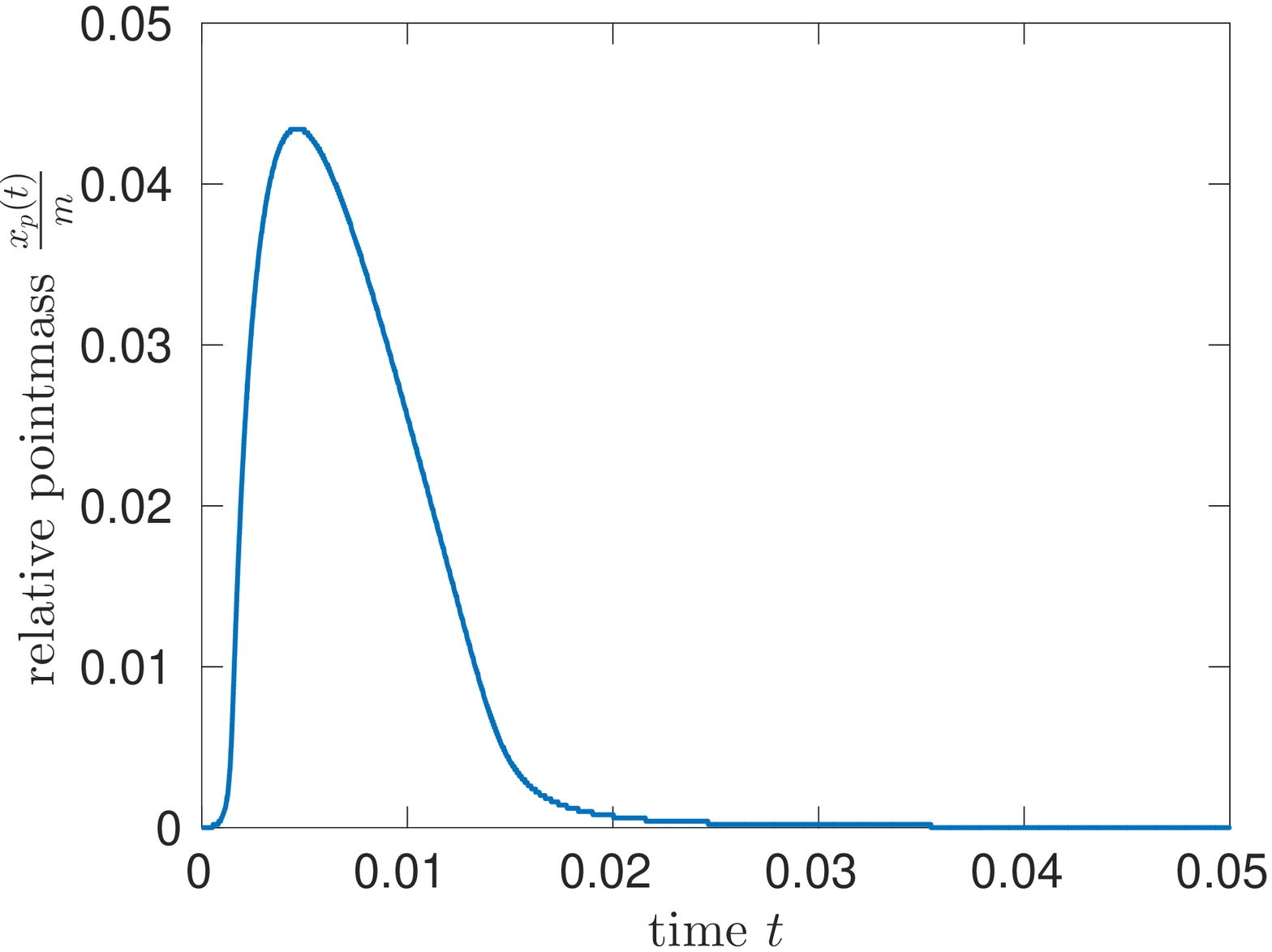} 
\subcaption{Dirac mass ($\sigma=0.1$).}
\label{fig:133}
\end{minipage}
\caption[Large-time behaviour for~\ref{it:sub} and~\ref{it:subbu} ($d=1, \gamma=2.9, m<m_c$)]{The mass-subcritical cases~\ref{it:sub} and~\ref{it:subbu}, $d=1, \gamma=2.9, A=1.5.$}
\label{fig:tc}
\end{figure}

\paragraph{Entropy decay.} The convergence to the minimiser of the entropy can be clearly observed in Figures~\ref{fig:111} and~\ref{fig:141}. Beyond,
Figures~\ref{fig:112},~\ref{fig:142},~\ref{fig:122} and \ref{fig:132}, which show the evolution of the relative entropy $H(u(t,\cdot))-H(u_{\infty})$, indicate an exponential decay of the entropy.
 The red slopes in Figures~\ref{fig:112},~\ref{fig:142},~\ref{fig:122} and \ref{fig:132} indicate the approximate slopes of the graphs averaged over the intervals where they are plotted. The computed slopes imply quantitative decays rates for the entropy of the form  $e^{-\alpha t}$ 
 with the following numerical values for~$\alpha$:
 $\alpha\approx23.7$ for~\ref{it:sup},
 $\alpha\approx23.8$ for~\ref{it:sub}, 
 $\alpha\approx23.1$ for~\ref{it:subbu}, 
and $\alpha\approx23.0$ for~\ref{it:supasy}. 
  \begin{remark}
   In the mass-subcritical case $m<m_c$ there exists $T=T(u_0)<\infty$ such that 
   the mapping $u(t,\cdot)$ has no critical point for $t>T$, so that the density $f(t,\cdot)$ of its inverse is smooth (see~\cite[Corollary~4.5]{carrillo_finite-time_2019}). In this case one can exploit the fact that the entropy functional of the bosonic Fokker--Planck equation in 1D coincides with that of a nonlinear diffusion equation with linear drift to which the theory developed in~\cite{unterreiter_entropy_2001} applies
   in order to deduce exponential decay of the entropy with rate $\alpha=2$ for $t\ge T$, i.e.
   \begin{align*}
     H(u(t,\cdot))-H(u_\infty)\le (H(u(0,\cdot))-H(u_\infty))\mathrm{e}^{-2t},\;t\ge t_0,\quad t\ge T.
   \end{align*}
 This idea was used before in~\cite{carrillo_1d_2008} for~1D~KQ. The rate of convergence in the general case is still open. 
 \end{remark}

\paragraph{Finite-time condensation for $m>m_c$.} The finite-time condensation in the mass-supercritical case is well confirmed by the simulations \ref{it:sup}\&\ref{it:supasy}. 
Recall that the condensate corresponds to the zero level set of $u(t,\cdot)$, which we numerically determine by the criterion $|u(t,\cdot)|<10^{-6}$.
Figure~\ref{fig:113} shows the time evolution of the condensed part relative to the (conserved) total mass. It clearly shows the onset of a condensate after some time $0<t\ll0.025$. 
Further figures depicting the formation of condensates are Fig.~\ref{fig:111},~\ref{fig:141} and~\ref{fig:143}. Interestingly, in Figure~\ref{fig:143} the fraction of mass in the condensate is not monotonic, illustrating that, even when above the critical mass, a previously formed condensate may partially dissolve.

\paragraph{Blow-up profile.} Figures~\ref{fig:115} and~\ref{fig:145} show the behaviour of $f(t,v)-f_c(v)$ for $0<v\ll R_1$ at the times $t=0.04$ and $t=0.1$. The figures indicate an error of the form
\begin{align}\label{eq:linearError}
  f(t,v)-f_c(v) = c_\pm(t)|v| + o(|v|)\quad\text{ as }v\to0\pm
\end{align}
for suitable constants $c_+(t), c_-(t)\in\mathbb{R}$, which, for asymmetric solutions, need not necessarily coincide. 
The asymptotic behaviour in~\eqref{eq:linearError} not only confirms the leading order spatial profile obtained rigorously in~\cite{carrillo_finite-time_2019} (see~\eqref{eq:profile}),
but indicates that the error with respect to $f_c$ may typically be of first order in 
$|v|$ and thus smaller than the order $1-2/\gamma$ ensured by formula~\eqref{eq:profile}. (A rigorous derivation of the improved error control can be found in~\cite{hopf_thesis}.)
Let us also mention that in both figures the solution $u(t,\cdot)$ is not uniformly close to $u_\infty$, so that the asymptotic behaviour of the density near the origin at the chosen times is not due to the fact that the long-time limit of the density equals~$f_c$.

\paragraph{Transient condensates.}

In Figure~\ref{fig:tc} the behaviour of a mass-subcritical, but initially very concentrated solution is compared to the solution emanating from a more spread out datum.
In both cases the entropy decays exponentially. Observe that in the case of high concentration, the solution forms a condensate in finite time which eventually vanishes again. We refer to this phenomenon as a \textit{transient condensate}. 
Recall that for $d=1$ and $\gamma>2$ the existence of transient condensates is known rigorously~\cite{carrillo_finite-time_2019}. The simulations based on~\ref{it:subbu} illustrate very explicitly how, after some finite time, the function $u(t,\cdot)$ begins to forms a flat part at the horizontal axis, which eventually disappears again as the solution converges to the smooth, non-degenerate equilibrium (cf.~Figure~\ref{fig:135}).

 \subsection{Validating KQ by means of explicit solutions in~2D}
 \label{ssec:val_2D}
 
 In the case $d=2$, KQ is $L^1$-critical and---as shown in~\cite{canizo_fokkerplanck_2016}---its isotropic form can be transformed in an explicit way to a linear Fokker--Planck equation, whose solutions are explicit by means of the fundamental solution for this problem in $\mathbb{R}^2$. Here we want to use these explicit solutions to validate the proposed numerical scheme. Since all simulations are performed on a finite domain with zero flux boundary condition, the solutions to KQ obtained upon this transformation are only approximations of the exact solutions to our problem.  However, we obtain a good approximation of the solutions in $B(0,R_1)\subset\mathbb{R}^2$ with zero flux provided $R_1$ is chosen sufficiently large. This is due to the fact that the exact solutions in $\mathbb{R}^2$ emanating from the chosen initial data (Gaussians) have exponential decay in $|v|$. The same is true for their derivative with respect to $v$, implying that on the boundary $\partial B(0,R_1)$ of a centred ball of large enough radius $R_1\gg1$ the flux is negligible. Hence, the exact solutions on $\mathbb{R}^2$ restricted to $B(0,R_1)$ are close to the exact solutions on $B(0,R_1)$ with zero flux.
 
 Next, we recall the transformation leading to the explicit formula of solutions on the whole space, as observed in~\cite{canizo_fokkerplanck_2016}:
 the solutions of the linear Fokker--Planck equation
 \begin{align}\label{eq:linFP}
   \partial_th & = \Delta h+\mathrm{div}(vh),\quad t>0,v\in\mathbb{R}^2,
   \\h(0,\cdot) & = h_0\nonumber 
 \end{align}
  are given by means of the fundamental solution 
  \begin{align*}
    F(t,v,w)=a(t)^{-1}K_{b(t)}(a(t)^{-1/2}v-w),
 \end{align*}
 where
 $ a(t)=\mathrm{e}^{-2t},\;\;b(t)=\mathrm{e}^{2t}-1,\text{ and }K_b(z)=(2\pi b)^{-1}\mathrm{e}^{-|z|^2/2b}.$
 More precisely, (for sufficiently regular data $h_0$) the solution of equation~\eqref{eq:linFP} takes the form
\begin{align}\label{eq:solFP}
  h(t,v) = \int_{\mathbb{R}^2} F(t,v,w)h_0(w)\,\d w.
\end{align}

 The relation between non-negative, isotropic solutions $f$ of $2$D KQ and non-negative, isotropic solutions $h$ of eq.~\eqref{eq:linFP} is given by 
 \begin{align}\label{eq:trafo}
   f(t,v) = \frac{h(t,v)}{1+\bar M_h(t,|v|)}\quad \text{ resp. }\quad h(t,v)=f(t,v)\mathrm{e}^{\bar M_f(t,|v|)},
 \end{align}
where 
\begin{align*}
  \bar M_f(t,\rho)=\frac{1}{2\pi}\int_{\{|v|\le \rho\}}f(t,w)\,\d w
  = \int_0^\rho g(t,r)r\,\d r.
\end{align*}
 
 We initialise our tests again with a centred Gaussian of the form
 \begin{align*}
   f_0(v) = A\mathrm{e}^{-\frac{|v|^2}{2\sigma^2}}
 \end{align*}
 for fixed positive constants $A$ and $\sigma$. Then the initial datum $h_0$ corresponding to $f_0$ via the transformation~\eqref{eq:trafo} is given by 
 \begin{align*}
   h_0(v) = A\mathrm{e}^{-\frac{|v|^2}{2\sigma^2}}\mathrm{e}^{A\sigma^2\left(1-\mathrm{e}^{-\frac{|v|^2}{2\sigma^2}}\right)},
 \end{align*}
 and from formula~\eqref{eq:solFP} and relation~\eqref{eq:trafo} we infer an expression for the solution $f$, which shows, in particular, that $f(T,\cdot)$ has exponential decay for any positive time $T$.  
 In our actual code, we use the inverse cdf of $f(T,\cdot)$.

 \paragraph{Details on the tests.}
 We choose $R_1$ to be the smallest radius satisfying $f_c(v)\le10^{-4}$ for $|v|\ge R_1$. This guarantees that for any not too large $\sigma>0$, the function $f(t,\cdot)$ is small outside $B(0,R_1)$.
  
 Two different tests are performed using the following common set of parameters: 
 $A=4$, $\sigma=0.9$, final time $T=0.04$ and size of the coarsest mesh equal to $n_0=25$.
 Since the solution to the exact problem remains bounded,  
the tests are performed with $\varepsilon=\delta=0$. 
 
 In the first test the dependence of the $L^2$ distance at time $T$ between exact and computed solution for different spatial resolutions is analysed. 
 More precisely, for $j=0,\dots,N=5$ we compute the error 
 $$E_j= \|S^{(j)}(T,\cdot)-S^{(j)}_\mathrm{exact}(T,\cdot)\|_{l^2(J_j)}\cdot 2^{-j},$$ where
 $J_j$ denotes the discrete mesh using a total number of $2^{j}n_0+1$ mesh points intersected with the interval $[0,m/2]$, $S^{(j)}_\mathrm{exact}$ denotes the exact solution restricted to the  spatial mesh $J_j$ and $S^{(j)}$ the discrete solution computed on the mesh $J_j$ using a total number of $400$ time steps.
 Since we expect a polynomial dependence of the error on the spatial increment, we then let $\mathrm{rate}(j)=\log_2(E_j/E_{j+1})$. 
 The results of the test can be found in Table~\ref{table:2Dmesh}. Theoretically, since in the present case of two space dimensions the original density $f$ remains uniformly bounded in time, which implies that $\partial_zS$ stays away from zero, the spatial discretisation based on central differences should guarantee a quadratic dependence of the truncation error on the spatial increment. The rates displayed in Table~\ref{table:2Dmesh} are somewhat worse, possibly due to the fact that the mesh size has not been chosen sufficiently large to capture the asymptotic behaviour well enough.
 
In the second test we analyse the dependence of the $L^2$ space-time distance between exact and computed solution on the number of spatial and temporal grid points. The procedure is analogous to the first test except that the $j$-th mesh is obtained by using $2^{j}n_0+1$ spatial and $2^{j}m_0$ temporal grid points, where $m_0=4$, and that now the error is given by
$$E_j= \|S^{(j)}-S^{(j)}_\mathrm{exact}\|_{l^2(I_j\times J_j)}\cdot 2^{-2j},$$
where $I_j$ denotes the discrete temporal mesh consisting of $2^{j}m_0$ time points. The results are displayed in Table~\ref{table:2Dfull} and suggest a linear rate of convergence. This is in line with the backward Euler scheme used for the time stepping.

\begin{table}[H]\footnotesize
\parbox{.5\linewidth}{
\centering
\begin{tabular}{|cccc|}                
\hline number of&mesh size&$L^2$ error& rate\\                                  
 time points &  & (at time $T$) &  \\   
\hline  \hline                                  
4000 & 25 & 6.2783e-3   & - \\ 
     \hline 
     4000 & 50 & 2.2323e-3  & 1.4919 \\ 
     \hline 
     4000 & 100 & 7.9661e-4 & 1.4866 \\ 
     \hline 
     4000 & 200 & 2.6080e-4 & 1.6109 \\ 
     \hline 
     4000 & 400 & 7.7921e-5 & 1.7428 \\ 
     \hline 
     4000 & 800 & 1.9283e-5 & 2.0147 \\ 
     \hline                                        
\end{tabular}    
     \caption{Convergence to exact solution at the final time $T=0.04$.}  
      \label{table:2Dmesh} 
}
\hfill
\parbox{.5\linewidth}{
\centering

      \begin{tabular}{|cccc|} 
    \hline     
    number of&mesh size& full $L^2$ error& rate\\                                  
 time points &  &  &  \\   
\hline  \hline  
  4 & 25 & 8.3850e-4 & - \\
    \hline
    8 & 50 & 4.1295e-4 & 1.0218 \\
    \hline
    16 & 100 & 2.0813e-4 & 0.9885 \\
    \hline
    32 & 200 & 1.0427e-4 & 0.9971 \\
    \hline
    64 & 400 & 5.1996e-5 & 1.0039 \\
    \hline
    128 & 800 & 2.5774e-5 & 1.0125 \\

      \hline  
      \end{tabular}  
      \caption{Convergence to reference solution (spacetime grid).}  
      \label{table:2Dfull} 
}
\end{table}

\begin{remark}[Validation of regularisation]
 For completeness, we also tested the dependence of the computed solution on the regularisation parameters $\varepsilon$ and $\delta$, even though this is not necessary for $2$D KQ since the density is theoretically known to remain bounded.
 We obtained a polynomial decrease of the error.
\end{remark}

\subsection{Simulations of 3D KQ in radial coordinates}\label{ssec:sim3D}

\begin{figure}[H]
\begin{minipage}{0.49\textwidth}\centering
\includegraphics[scale=\sca]{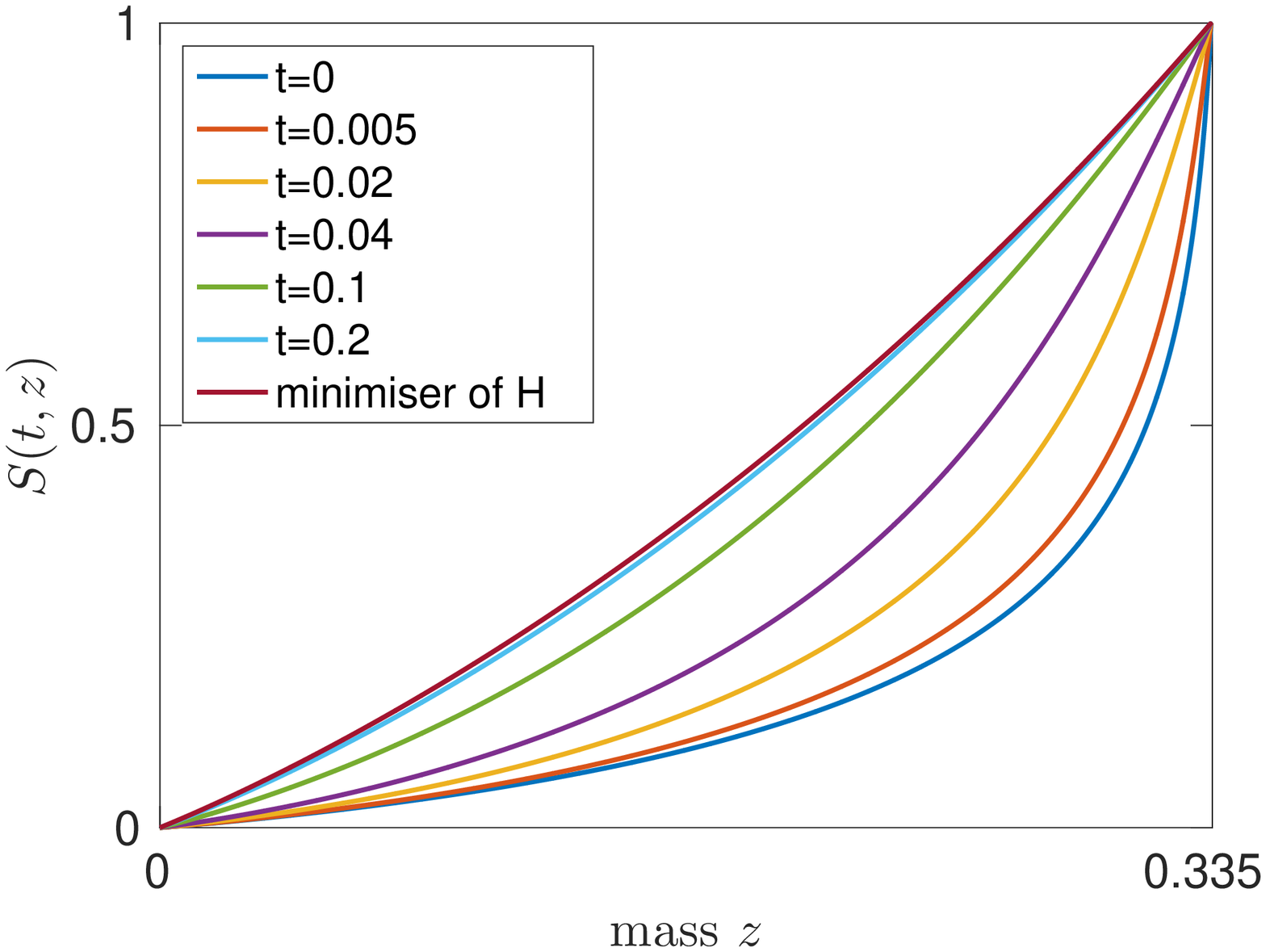}  
\subcaption{$S(t,\cdot)$ and $S_\infty$.}
\label{fig:321}\vfill
\end{minipage}
\begin{minipage}{0.49\textwidth}\centering
  \includegraphics[scale=\sca]{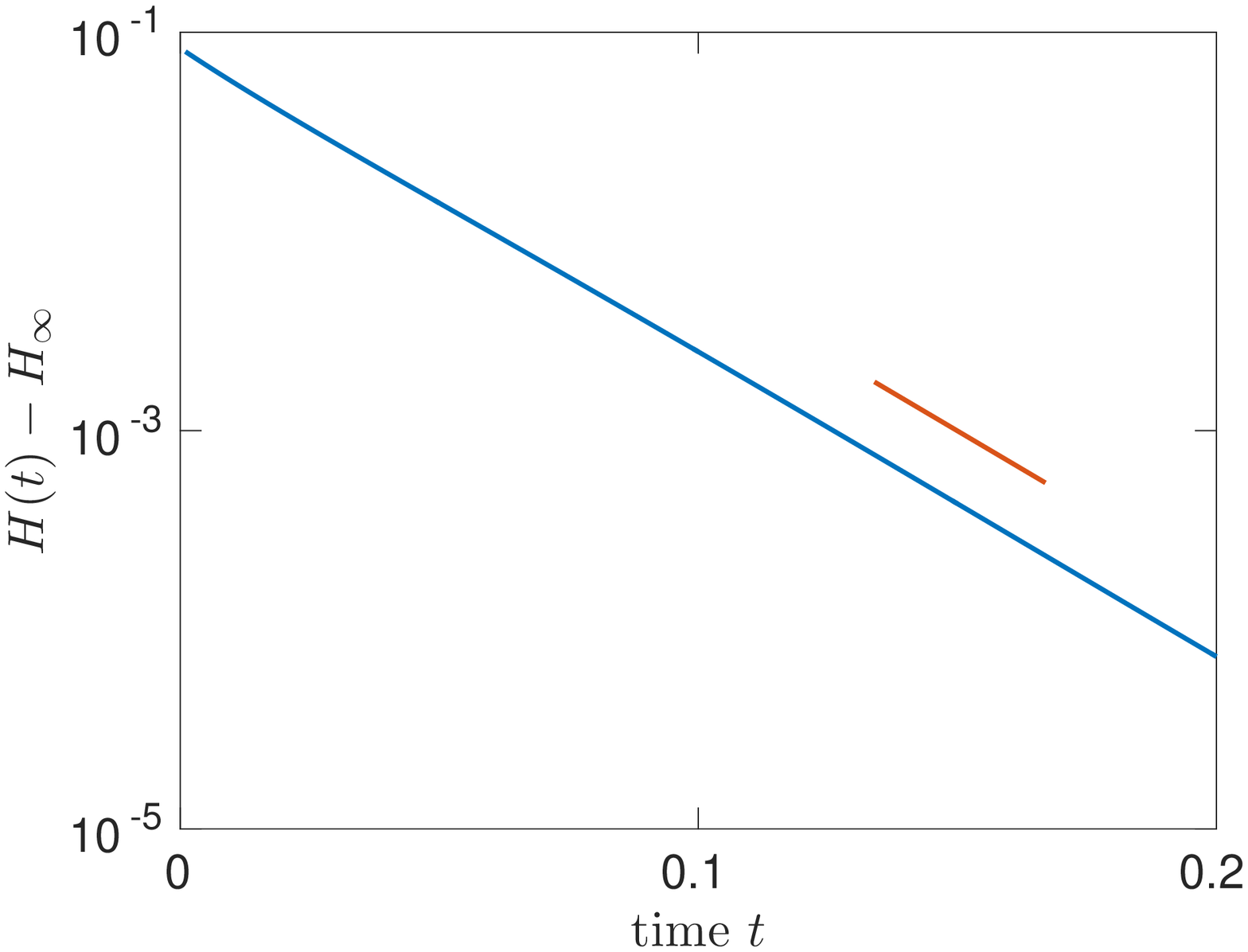}
  \subcaption{Evolution of the relative entropy.}
\label{fig:322}
\end{minipage}
\caption[Large-time behaviour for~\ref{it:sub3D} ($d=3, \gamma=1, m<m_c$)]{Long-time behaviour in mass-subcritical case~\ref{it:sub3D} ($\gamma=1, d=3$).}
\end{figure}

\begin{figure}[H]
 
\begin{minipage}{0.49\textwidth}\centering
\includegraphics[scale=\sca]{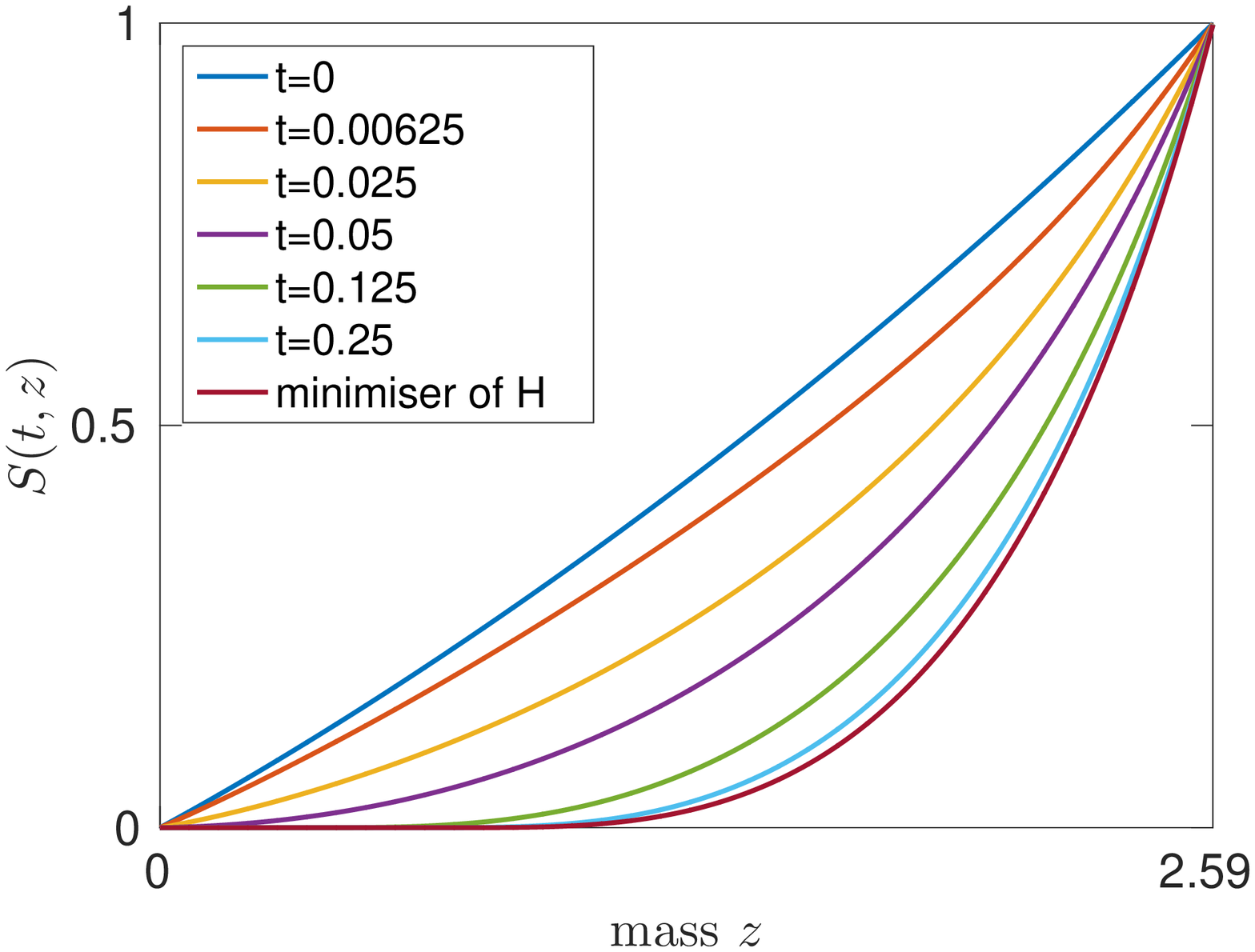} 
\subcaption{$S(t,\cdot)$ and $S_\infty$.}
\label{fig:311}
\end{minipage}
\begin{minipage}{0.49\textwidth}\centering
\includegraphics[scale=\sca]{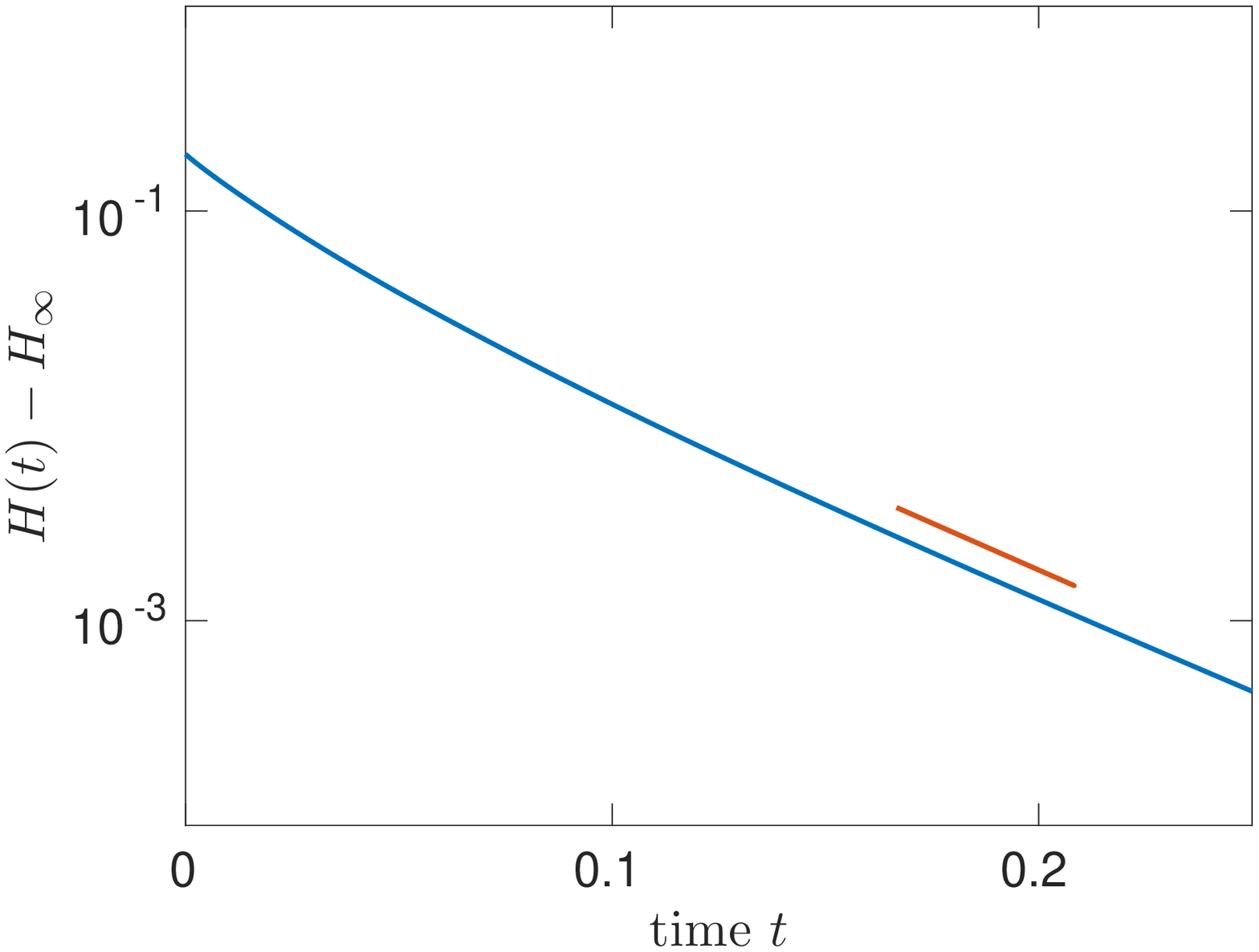}
\subcaption{Evolution of the relative entropy.
}
\label{fig:312}
\end{minipage}
\vspace{0.5cm}

\begin{minipage}{0.49\textwidth}\centering
\includegraphics[scale=\sca]{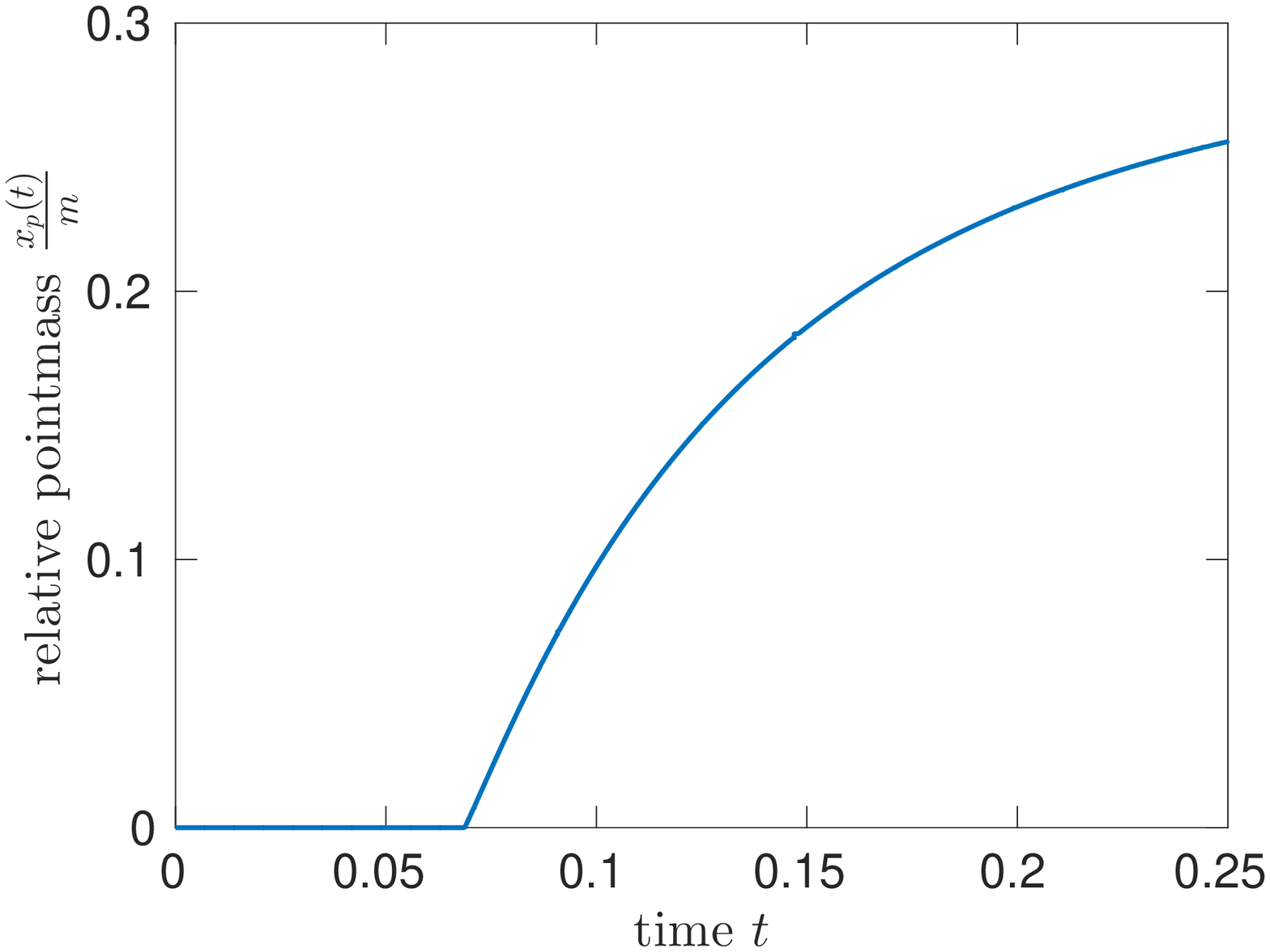} 
\subcaption{Evolution of the Dirac part.}
\label{fig:313}
\end{minipage}
\begin{minipage}{0.49\textwidth}\centering
\includegraphics[scale=\sca]{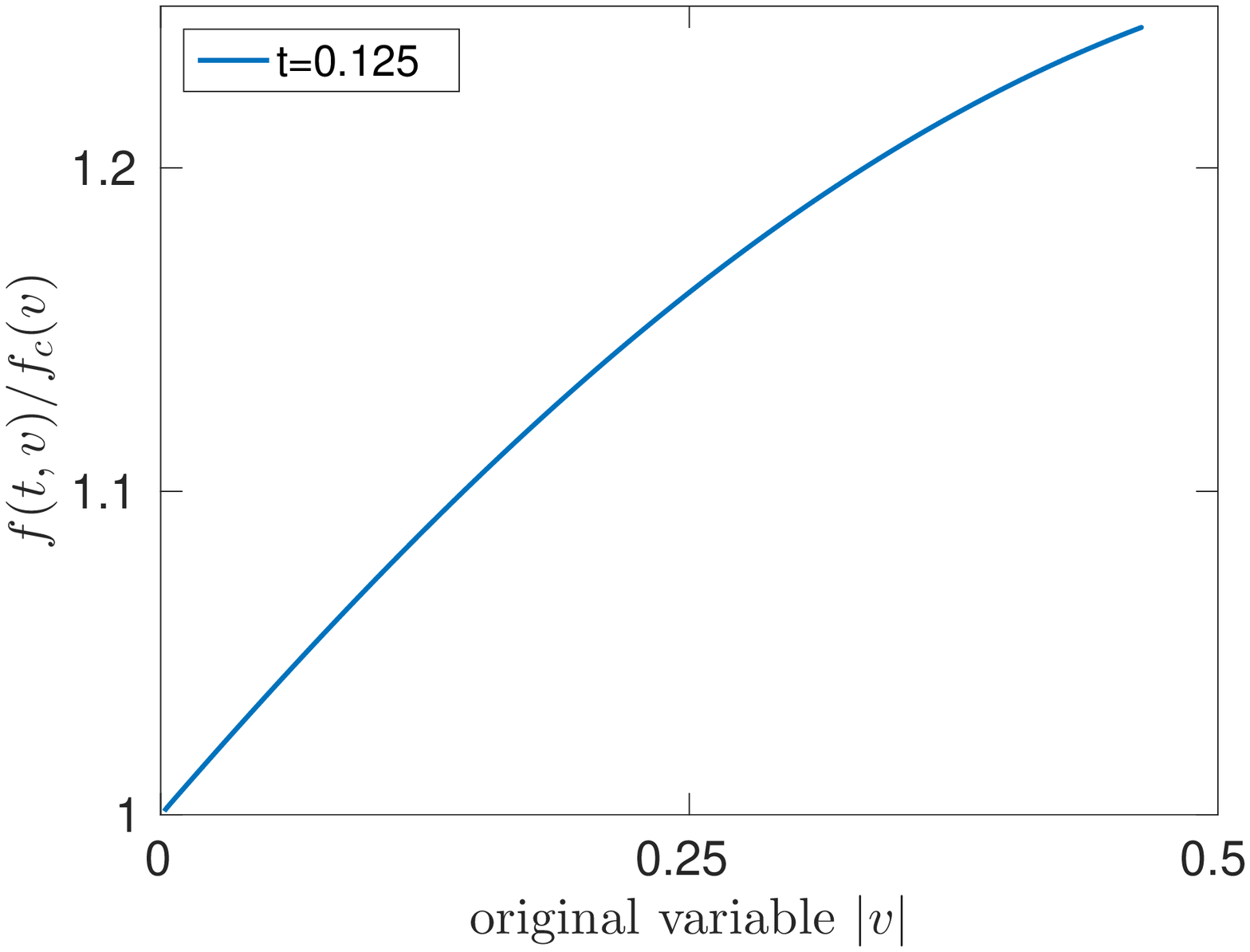} 
\subcaption{Behaviour near singularity.
}
\label{fig:314}
\end{minipage}
\caption[Large-time behaviour for~\ref{it:sup3D} ($d=3, \gamma=1, m>m_c$)]{Long-time behaviour in the mass-supercritical case~\ref{it:sup3D} ($d=3, \gamma=1, \varepsilon=10^{-12}$ and $ \delta=0$).}
\end{figure}

\begin{figure}[H]
  \begin{minipage}{0.49\textwidth}\centering
\includegraphics[scale=\sca]{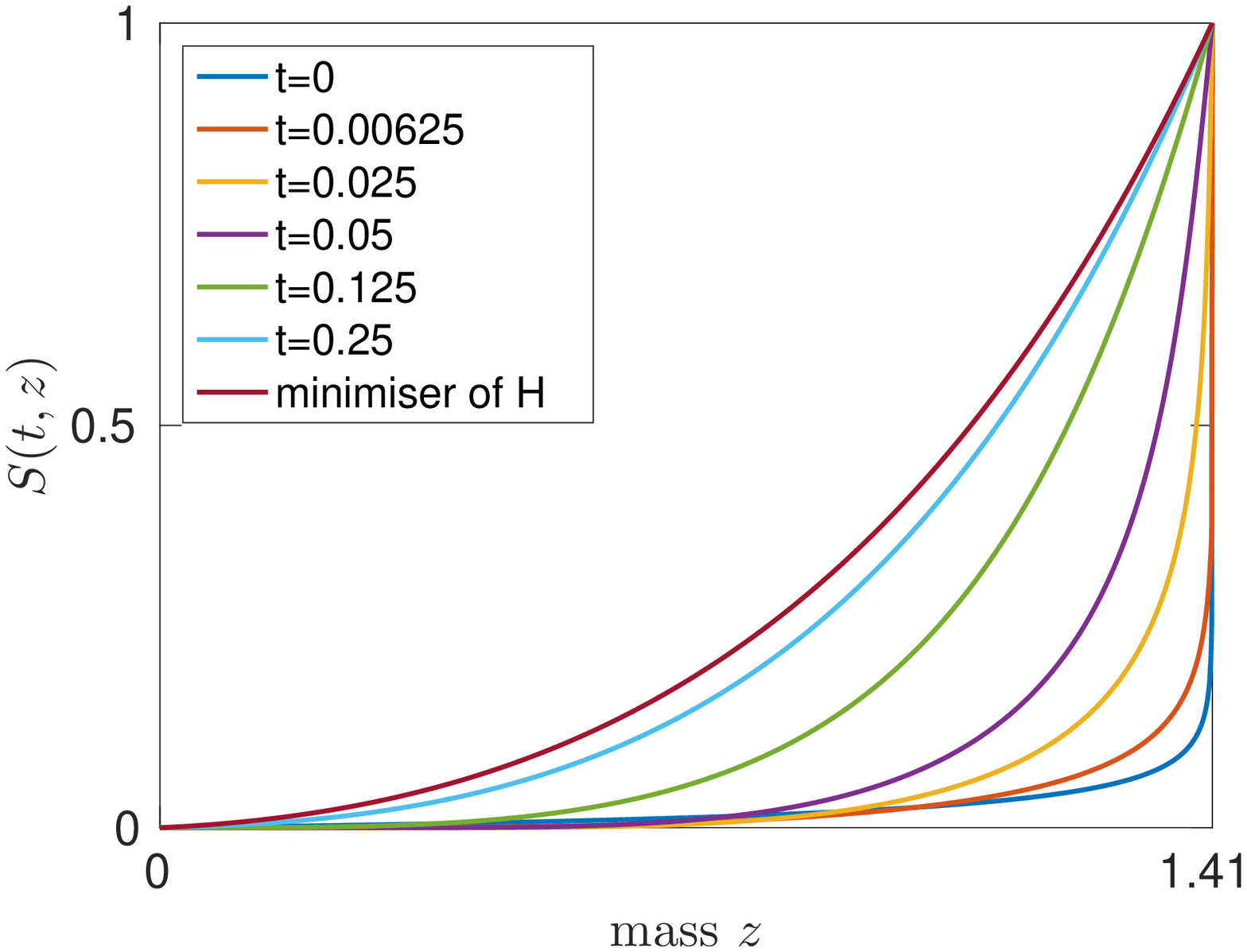}
\subcaption{$S(t,\cdot)$ and $S_\infty$.}
\label{fig:381}
\end{minipage}
  \begin{minipage}{0.49\textwidth}\centering
\includegraphics[scale=\sca]{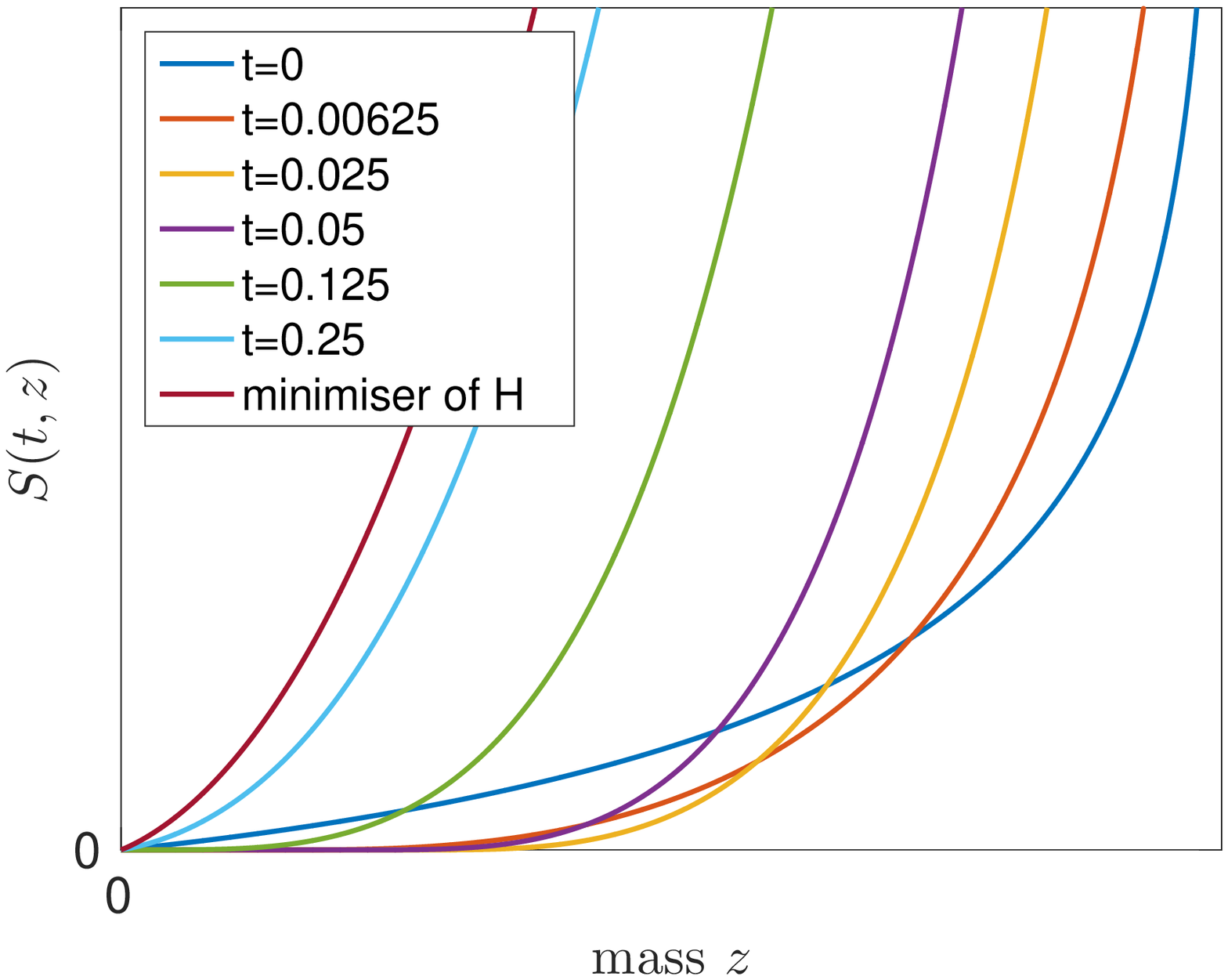}
\subcaption{Zoomed-in view of Fig.~\ref{fig:381}.}
\label{fig:385}
\end{minipage}

\vspace{.5cm}

  \begin{minipage}{0.49\textwidth}\centering
\includegraphics[scale=\sca]{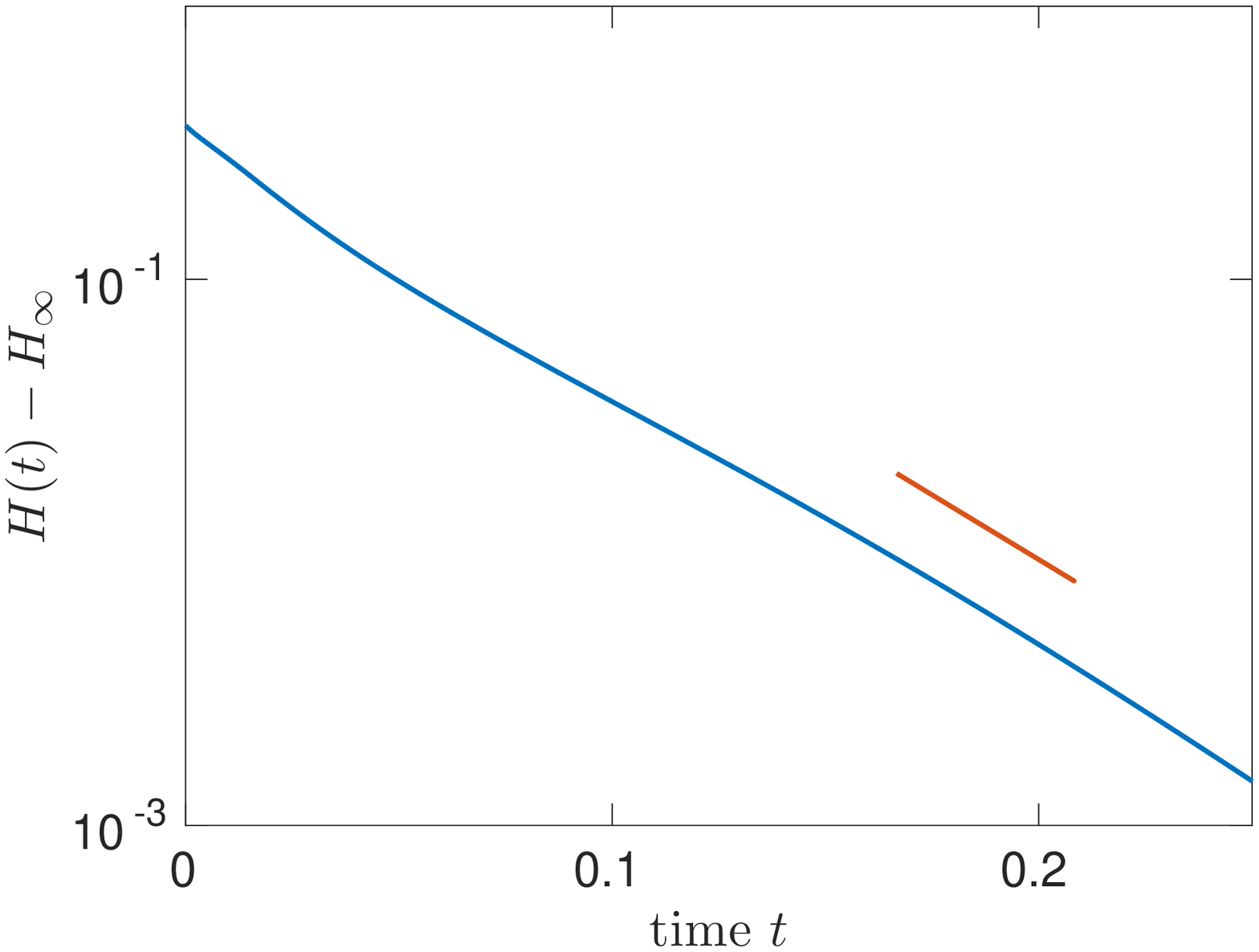}
\subcaption{Evolution of the relative entropy.
}
\label{fig:382}
\end{minipage}
 \begin{minipage}{0.49\textwidth}\centering
\includegraphics[scale=\sca]{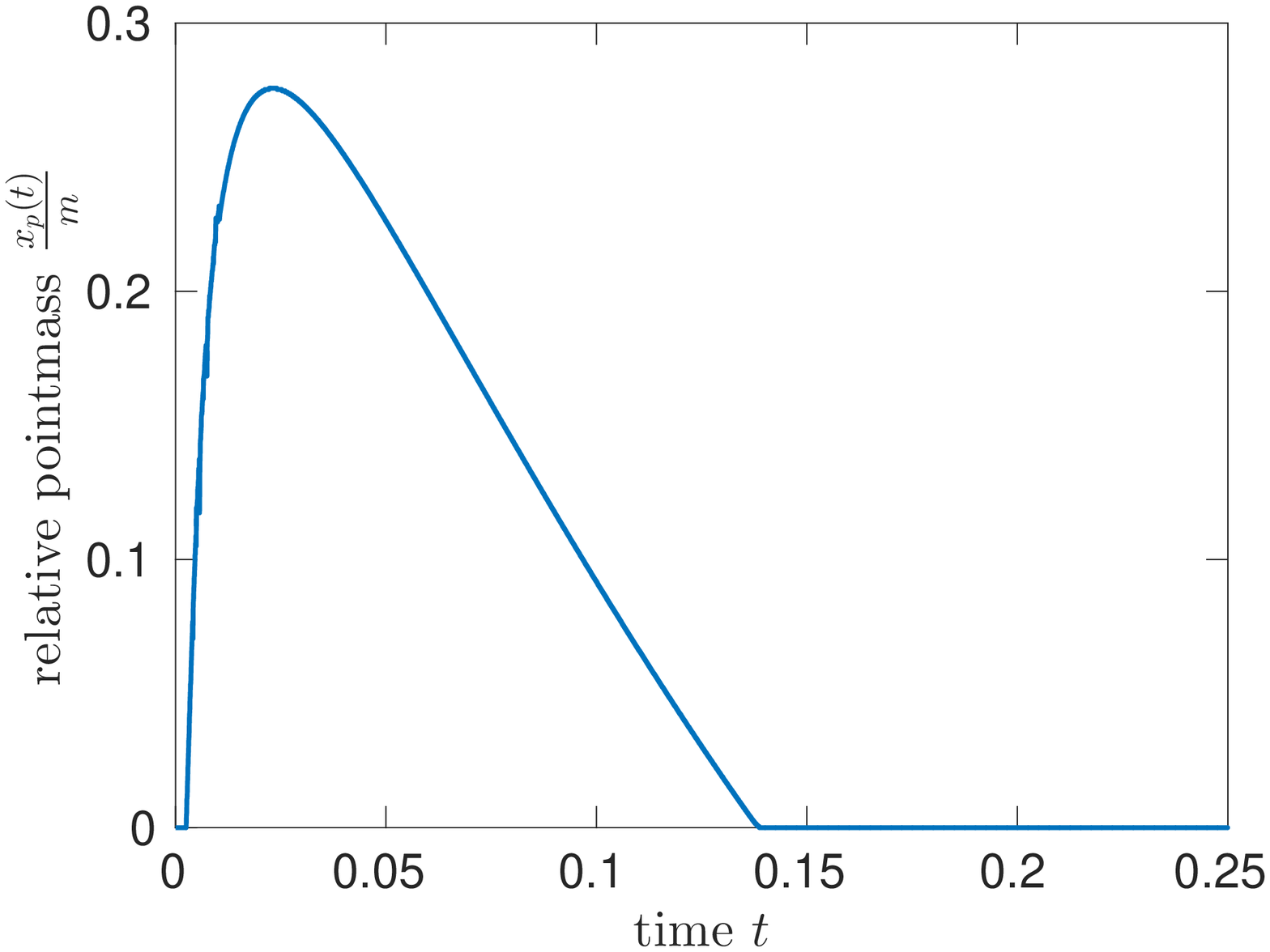}
\subcaption{Evolution of the Dirac part.}
\label{fig:383}
\end{minipage}
\caption[Transient condensate for~\ref{it:subbu3D} ($d=3, \gamma=1, m<m_c$)]{Transient condensate in the mass-subcritical case~\ref{it:subbu3D} ($d=3, \gamma=1, \varepsilon=\delta=10^{-10}$).}
\label{fig:tc3d}
\end{figure}

Here, we simulate equation~\eqref{eq:2regS3D} with $d=3$ for suitable choices of $\varepsilon,\delta$, $0<\varepsilon,\delta\ll1$, where we choose $R_1=1$. 
We recall our notation $\bar m_c=\frac{1}{|\partial B(0,1)|}\int_{B(0,R_1)}f_c(v)\,\d v$, where now $|\partial B(0,1)| = 4\pi$ denotes the area of the 2-sphere, and remark that the numerical value of $\bar m_c$ is approximately given by~$\bar m_c\approx1.84$.
We perform three simulations with a mass-supercritical, a mass-subcritical and a highly concentrated initial datum, respectively. More precisely, choosing as initial data again Gaussians of the form $f_0(v)=A\mathrm{e}^{-|v|^2/(2\sigma)}$, we run our scheme with the following three sets of parameters:
\setlist[enumerate,1]{start=5}
\begin{enumerate}[label=(P\arabic*)]   
  \item\label{it:sub3D} $m<m_c:$  $\sigma=0.3$, $A=3$, $T=0.2$, $\tau=0.001$, $n=2001$, $\varepsilon=0$, $\delta=0$.
  \item\label{it:sup3D} $m>m_c:$ $\sigma=0.9$, $A=10$, $T=0.25$, $\tau=5\cdot10^{-6}$, $n=50001$,  $\varepsilon=10^{-12}$, $\delta=0$.
 \item\label{it:subbu3D} $m<m_c:$  $\sigma=0.15$, $A=50$, $T=0.25$, $\tau=5\cdot10^{-5}$, $n=2001$, $\varepsilon=10^{-10}$, $\delta=10^{-10}$.
\end{enumerate}
The quantity $\bar m:=m/|\partial B(0,1)|$ associated with the above choice of parameters takes the value $\bar m\approx 0.335$ for~\ref{it:sub3D}, $\bar m\approx 2.59$ for~\ref{it:sup3D}, and $\bar m\approx 1.41$ for~\ref{it:subbu3D} (see Figures~\ref{fig:321}, \ref{fig:311} and \ref{fig:381}).

The size of the condensate divided by $|\partial B(0,1)|$, i.e.~$\bar x_p(t):=\mathcal{L}^1(\{S(t,\cdot)=0\})$, is numerically determined by replacing the condition $S(t,\cdot)=0$ with the smallness criterion $S(t,\cdot)<10^{-10}$.
\begin{remark}
The choice of the comparatively fine mesh in~\ref{it:sup3D} was made in order to ensure a sufficiently good approximation of the evolution of the entropy. See Fig.~\ref{fig:312}, which suggests an exponential decay.
\end{remark}

\paragraph{Long-time behaviour.}
Our simulations suggest that $3$D KQ has properties which are very similar to the Fokker--Planck model for bosons in $1$D in the $L^1$-supercritical regime.  Figures~\ref{fig:321},~\ref{fig:311} and~\ref{fig:381} suggest that in the long-time limit the numerical solution $S(t,\cdot)$ approximates the minimiser of the entropy (at the level of $S$), which we here denote\footnote{For simplicity, in our notation $S_\infty$ for the entropy minimiser we omit its dependence on the given mass $m$ and the radius $R_1$. } by $S_\infty$.
\paragraph{Entropy.}
The decay of the relative entropy appears to be exponential in all three cases~\ref{it:sub3D}--\ref{it:subbu3D}, see Figures~\ref{fig:322},~\ref{fig:312} and~\ref{fig:382}. In each of these plots the red slope indicates the approximate slope of the graph averaged over the interval where it is plotted. Numerically, the relative entropy $H(t)-H_{\infty}$ appears to decay to zero like $e^{-\alpha t}$, where 
 $\alpha\approx35.3$ for~\ref{it:sub3D}, $\alpha\approx21.1$ for~\ref{it:sup3D}, and $\alpha\approx21.7$ for~\ref{it:subbu3D}.

\paragraph{Condensation.}
In both the mass-supercritical case~\ref{it:sup3D} and the case of high concentration near the origin~\ref{it:subbu3D} we observe the onset of a flat part at the level of  $S(t,\cdot)$ at height zero after some finite time, see Fig.~\ref{fig:313} and~\ref{fig:383}. In the original variables this means that mass is gradually absorbed by the origin. Furthermore, Fig.~\ref{fig:383} shows that, similarly to the observations in $1$D (see Section~\ref{ssec:sim1D}), it is possible for mass previously concentrated at velocity zero to escape. In fact, the condensate component may even dissolve completely. Thus, at least numerically, the fraction of particles in the condensate is, in general, not monotonic in time for the 3D Kaniadakis--Quarati model.

\begin{remark}
	In order to produce the transient condensate in Figure~\ref{fig:tc3d} it was necessary to choose the parameter $\delta$ appearing in equation~\eqref{eq:2regS3D} (and its discrete counterpart) 
	strictly positive.
	The same simulation for $\delta=0$ results in the flat part being trapped at height zero once it has formed. As explained in Section~\ref{sssec:eqhd} and also in view of our results for the $1$D model, this \enquote{stickiness} appears to be a numerical artefact resulting from the circumstance that a regularisation based on a positive $\varepsilon$ but vanishing~$\delta$ is imbalanced and favours condensation.
\end{remark}

\begin{figure}[H]\centering
\begin{minipage}{0.49\textwidth}\centering
\includegraphics[scale=\sca]{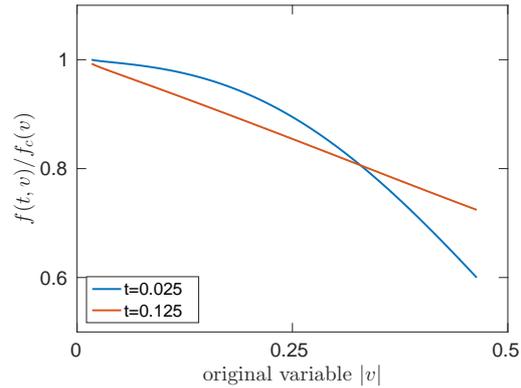} 
\vfill
\end{minipage}
\caption[Blow-up profile in \ref{it:subbu3D} ($d=3, \gamma=1,m<m_c$)]{Spatial blow-up profile in \ref{it:subbu3D}.}
\label{fig:384}
\end{figure}

\paragraph{Blow-up profile.} At times where the solution has a non-trivial condensate component, we were interested in the spatial behaviour of $S(t,\cdot)$ close to $\{S(t,\cdot)=0\}$. 
Owing to the results on the $1$D model, one may expect the function $f(t,\cdot)$ to behave to leading order like the limiting steady state $f_c$, i.e.\;like $2|v|^{-2}$. 
Furthermore, the formal expansions in~\cite[Section~III.C]{sopik_dynamics_2006} suggest that for isotropic solutions of 3D KQ the error by which $f(t,\cdot)$ deviates from $f_c$ has the form
\begin{align}\label{eq:Sopik3Dprofile}
  f(t,v) - f_c(v) = c(t)|v|^{-1} + o(|v|^{-1})
\end{align}
for some constant $c(t)\in\mathbb{R}$.
Our experiments corroborate formula~\eqref{eq:Sopik3Dprofile}. Indeed, Figures~\ref{fig:314} and~\ref{fig:384} displaying the quantity $f(t,v)/f_c(v)$ at times where $f(t,\cdot)$ is unbounded at the origin show that numerically it behaves like $1+\tilde c(t)|v|+o(|v|)$ as $|v|\to0$.

\section{Conclusion}
In this work we propose a numerical scheme for nonlinear Fokker--Planck equations for bosons able for the first time to cope with Dirac delta concentrations of (partial) mass at the origin in finite time and to go beyond this blow-up time. This is achieved by
considering appropriately normalised pseudo-inverse distributions and scaling suitably the equation to obtain an alternative formulation admitting a Dirac delta concentration at the origin as a possible steady state. These new PDEs are solved by implicit schemes, and their approximations by Newton--Raphson type methods are shown to be numerically convergent by mesh refinement, even beyond the blow-up time. 
The physical entropy associated to these problems is shown to be decreasing for the semidiscrete schemes in 1D and 2D. 
We illustrate different phenomena appearing in the 3D radial KQ model mimicking the phenomena observed and partially proved for the 1D caricature of the KQ model in the $L^1$-supercritical case, see \cite{carrillo_finite-time_2019}.

\section*{Acknowledgements}
JAC was partially supported by the EPSRC grant number EP/P031587/1. KH was supported by MASDOC DTC at the University of Warwick, which is funded by the EPSRC grant number EP/HO23364/1. MTW acknowledges partial support by the EPSRC grant number EP/P01240X/1.

\bibliographystyle{abbrv}
\bibliography{bib/sim}

\end{document}